\documentclass[11pt]{amsart}
\usepackage{amscd,amsthm,amsfonts,amssymb,esint,tikz-cd}
\usepackage{mathtools}
\usepackage[colorlinks=true]{hyperref}
\usepackage{fullpage}
\usepackage[all]{xy}
\usepackage{graphicx}

\usepackage{hyperref}
\newtheorem{theorem}{Theorem}[section]
\newtheorem{corollary}[theorem]{Corollary}
\newtheorem{lemma}[theorem]{Lemma}

\newtheorem{proposition}[theorem]{Proposition}
\newtheorem{question}[theorem]{Question}
\newtheorem{remark}[theorem]{Remark}

\newcommand{\ZZ}{{\mathbb Z}}

\newcommand{\tr}{{\rm tr}}

\newcommand{\id}{{\rm id}}

\newcommand{\degr}{{\rm deg}}
\newcommand{\MCG}{{\rm MCG}}
\newcommand{\MSO}{{\rm MSO}}
\newcommand{\HNN}{{\rm HNN}}
\newcommand{\Tor}{\rm{Tor}}

\begin{document}

\title[]
{$\pi_1$-injective bounding \\ and application to 3- and 4-manifolds}

\author{Jianfeng Lin}
\address{Department of Mathematics Science, Tsinghua University, Beijing, 100080, CHINA}
\email{linjian5477@mail.tsinghua.edu.cn}

\author{Zhongzi Wang}
\address{Department of Mathematical Sciences, Peking University, Beijing 100871 CHINA}
\email{wangzz22@stu.pku.edu.cn}

%\author{......}



\thanks{The first author is partially supported by Simons Collaboration Grant 615229.}

\begin{abstract} 
Suppose a closed oriented $n$-manifold $M$ bounds an oriented $(n+1)$-manifold.
It is known that  $M$ $\pi_1$-injectively bounds an oriented $(n+1)$-manifold  $W$.
We prove that   $\pi_1(W)$ can be  residually finite if  $\pi_1(M)$  is,
and  $\pi_1(W)$ can be finite if  $\pi_1(M)$ is. In particular,
 each closed 3-manifold $M$ $\pi_1$-injectively bounds a 4-manifold
with residually finite $\pi_1$, and  bounds a 4-manifold with finite $\pi_1$ if $\pi_1(M)$ is finite.
Applications to 3- and 4-manifolds are given:

(1) We study finite group actions on closed  4-manifolds 
and $\pi_1$-isomorphic cobordism of 3-dimensional lens spaces. Results including:
(a) Two lens spaces are $\pi_1$-isomorphic cobordant if and only if there is a degree one map between them.
(b) Each spherical 3-manifold $M\ne S^3$ can be realized as the unique non-free orbit type for  a finite group action on a closed 4-manifold. 

(2)  The minimal bounding index $O_b(M)$ for closed 3-manifolds $M$ are defined, %and bounding Euler charicteristic  $\chi_b(M)$.
the relations between  finiteness of  $O_b(M)$ and virtual achirality of aspherical (hyperbolic) $M$ are addressed.
We calculate  $O_b(M)$ for some lens spaces $M$. Each prime is realized as a minimal bounding index.

(3) We also discuss some concrete examples:
Surface bundle often bound surface bundles, and prime 3-manifolds often
 virtually bound surface bundles, $W$ bounded by some lens spaces realizing $O_b$ is constructed.
\end{abstract}

%\vskip 1true cm

%{\bf Key words:} 3-manifolds; 4-manifolds;  fundamental groups.

\date{} 
\maketitle

\tableofcontents

\section{Introduction}

% all finite group actions on manifolds are 
%smooth， faithful, unless otherwise stated.

%The references for 3-manifold topology see \cite{He}, for 3-manifold geometry see \cite{Sc}, for algebraic topology see \cite{Ha}.

%\subsection{Manifolds with finite $\pi_1$ bound manifolds with finite $\pi_1$}
All manifolds discussed in this paper are oriented and compact.
Suppose $M_1$ and $M_2$ are closed oriented $n$-manifolds.
Say $M_1$ and $M_2$  are cobordant, if there is an oriented $(n+1)$-manifold $W^{n+1}$ such that 
 
$$\partial W=M_1\coprod -M_2$$
where $\partial W$ is the boundary of $W$ and $-M$ is the manifold $M$ but with an opposite 
orientation. The cobordant relation of closed oriented $n$-manifolds is an equivalence relation,
and the equivalence classes form an abelian group under  disjoint  union, denoted by 
$\Omega_n$. 
The study  of $\Omega_n$ and  its various extensions is an important topic in topology with long history and is still attractive today, 
see \cite{Roh}, \cite{Tho}, \cite{CG} and \cite{DHST} for a few examples.
Call a closed oriented $n$-manifold $M$ bounding, if there is  compact oriented $(n+1)$-manifold $W$
such that $\partial W=M$.

%\begin{itemize}

%\item $\Omega_1=\Omega_2=0$ are classical facts.

%\item  Rohlin proved that $\Omega_3$ is trivial and $\Omega_4\cong \mathbb{Z}$ \cite{Roh}, for $n=3$, also \cite{Li}.

%\item Thom established  co-bordism theory \cite{Tho}, and many important studies follows after.

%\end{itemize}

%Each connected  CW-complex has a fundamental group $\pi_1(X)$.
Suppose $X_1$ and $X_2$ are connected CW-complexes.  Call a map $f: X_1\to X_2$ $\pi_1$-injective ($\pi_1$-surjective, $\pi_1$-isomorphic, respectively),
if the induced map  $f_*: \pi_1(X_1)\to \pi_1(X_2)$ is an injection (surjection, isomorphism, respectively). 
$\pi_1$-injective  embeddings of surfaces into 3-manifolds is a basic tools  in 3-manifolds, and 3-manifolds are almost determined by their $\pi_1$  \cite{He} \cite{Thu}.

%Fundamental group is fundamental in topology.

%\begin{itemize}

%\item Closed and orientable 1-man\citeifolds and  2-manifolds  are determined by their  $\pi_1$ are classical facts.

%\item Based on the theory of Thurston, closed and orientable 3-manifolds are almost determined by their $\pi_1$, and 
%3-manifold groups have some special properties, which  are important in the study of 3-manifolds  and often inspire the study of group theory.

%\item For each $n\ge 4$, a simple but important result is that any finitely presented group can be realized as  $\pi_1(M)$ for some $n$-manifold $M$.

%\end{itemize}
%\subsection{Motivations and Main results}

 It is known  that each closed oriented 3-manifold $M$ is bounding \cite{Roh} and 3-manifold
groups have many good properties \cite{Thu}. The following question is the main motivation of our study:

\begin{question}\label{motivation} Could $M$
bounds 4-manifold  $W$ $\pi_1$-injectively? Moreover, given some property $P$ of $M$, can we require that $W$ also has property $P$?
\end{question}
One can also ask this question in other dimensions. Indeed, some remarkable results in this issue have existed for a while:

%Before we state our results, we make some comments:

%One ask stronger question that could $M$ bounds a $W$ $\pi_1$-isomorphically? But this question often has no positive answer.

%begin{itemize}

%\item For $n=1$ the answer is NO due to the classification of 2-manifolds;

%\item For $n=2$ the answer is NO  by 3-manifold theory \cite[Theorem 10.2]{He};

%\item For $n=3$ the answer is NO unless $M$ is a connected sum of $S^2\times S^1$'s, \cite{Da}, also \cite{SW}.

%\item For $n\ge 4$ the answer is NO if some prime factor of $M$ is aspherical \cite{SW}.

%\end{itemize}

%We also make a comment on $\pi_1$-surjective bounding, even it gives no information about the $\pi_1$ of the bounded manifold.

%\begin{itemize}

%\item Each closed orientable $n$-manifold $M$, if its bounds,  it  $\pi_1$-surjectively bounds an orientable $(n+1)$-manifold,
%indeed it bounds a simply connected manifold if $n\ne 2$.  
%\end{itemize}

%Finally  we make comments on $\pi_1$-injective bounding:

\begin{itemize}

\item Hausmann proved that every closed oriented bounding $n$-manifold $\pi_1$-injectively bounds an  orientable $(n+1)$-manifold \cite{Hau} in 1981.

\item Davis-Januszkiewicz-Weinberger proved that every closed oriented aspherical bounding  $n$-manifold $\pi_1$-injectively bounds an   oriented aspherical $(n+1)$-manifold \cite{DJW} in 2001. 

\item Foozwell-Rubinstein proved that every closed Haken $3$-manifold $\pi_1$-injectively bounds a Haken $4$-manifold \cite{FR} in 2016.
\end{itemize}

\subsection{Statement of the main results} We are going to state our results for Question \ref{motivation}.
Call a group $G$  residually finite if for each $1\ne g\in G$, there is a finite group $H$ and a homomorphism $\phi: G\to H$
such that $\phi(g)\ne 1$. 
The  residually finite property is fundamental in the study of various virtual properties (we will discuss soon) of 3-manifolds  and of the theory of profinite groups.

\begin{theorem}\label{main}
Suppose $M$ is a closed oriented bounding $n$-manifold. Then 

{\rm (1)} $M$ $\pi_1$-injectively bounds a compact oriented $(n+1)$-manifold
with residually finite $\pi_1$ if $\pi_1(M)$ is residually finite.

{\rm (2)} $M$ $\pi_1$-injectively bounds a compact oriented $(n+1)$-manifold with finite $\pi_1$ if $\pi_1(M)$ is finite.
\end{theorem}

In order to prove Theorem \ref{main},  we will give an alternative proof  of Hausmann's Theorem \cite{Hau} in Section 3.

Since each closed 3-manifold has a residually finite $\pi_1$ \cite{Thu},  and $\Omega_3=0$ \cite{Roh}, we have 

\begin{theorem}\label{main3} Let  $M$ be a closed oriented $3$-manifold.

{\rm (1)}  $M$ $\pi_1$-injectively bounds a compact oriented $4$-manifold
with residually finite $\pi_1$. 

{\rm (2)} $M$ $\pi_1$-injectively bounds a compact oriented $4$-manifold
with finite $\pi_1$ if $\pi_1(M)$ is finite.
\end{theorem}

%{Applications to finite group actions on 4-manifolds}

Before discussing applications of Theorem \ref{main3} to 3-manifolds and 4-manifolds,
 we recall 
 Thurston's picture on 3-manifolds \cite{Thu}:
Let $Y$ be a  closed orientable prime 3-manifold. Then 
(i) $Y$ is either a $\mathcal G$-manifold,  where $\mathcal G$ is one of the following eight geometries:
$H^3$,  $Sol$, $Nil$,  $\widetilde{P
SL}(2,\mathbb{R})$, $S^3$, $E^3$,  $H^2\times E^1$, $S^2\times E^1$,
and $H^n, E^n, S^n$ indicate the $n$-dimensional hyperbolic, Euclidean, and spherical geometries;
or  (ii) $Y$ has a non-trivial JSJ tori decomposition such that each JSJ-piece of $Y$ supports the geometry of 
either  $H^2\times E^1$ or $H^3$,  and call $Y$  mixed if  at least one JSJ-piece is hyperbolic.

\subsubsection{Applications on group actions on $4$-manifolds}
The first application concerns finite group actions on 4-manifolds with prescribed orbit types. Suppose $G$ is a finite group acting on a closed, orientable 4-manifold $X$  whose non-free points are isolated.
Then we have the quotient map $q: X\to X/G$, and 
the $q$-image of each non-free orbit in $X/G$ has a neighborhood homeomorphic to a cone over a spherical 3-manifold $Y$. We call this $Y$ the type of this non-free orbit.
We say the $G$-action is semi-free, if it is free on the complement of its fixed points. 
And we say the $G$-action is almost free, if it has only one non-free orbit. One may ask the following natural questions.

\begin{question}\label{orbit type}
{\rm (1)} Which orbit types can arise from an almost free action on a $4$-manifold?

{\rm (2)} Which combinations of orbit types can arise from a semi-free action on a $4$-manifold?
\end{question}
Based on their fixed point theorem, Atiyah and Bott proved that two lens spaces are $h$-cobordant if and only if they are diffeomorphic \cite{AB}. One may wonder what happens if we weaken the condition to being $\pi_1$-isomorphic cobordant. Our theorem below answers this question. 

\begin{theorem}\label{appl-2}  Two lens spaces are $\pi_1$-isomorphic corbordant %,  or equivalently, there are types of a semifree $\mathbb{Z}/n$ action  on a simply connected closed 4-manifolds with only two fixed points, 
 if and only if there is an {orientation preserving homotopy equivalence between them}.
 \end{theorem}

Theorem \ref{appl-2}  follows from Theorem \ref{appl1} below, which answers the more general Question \ref{orbit type} (2) in the cyclic case. We use $\mathbb{Z}_n$ to denote the cyclic group of order $n$.

\begin{theorem}
\label{appl1} Let  $L(n, q_1), ... , L(n, q_m)$  be oriented lens spaces.  The following conditions are equivalent:
\begin{enumerate}
    \item There is a compact oriented $4$-manifold $W$ such that  $\partial W= \bigcup _{i=1}^m L(n, q_i)$ and each inclusion $L(n, q_i)\to W$ is $\pi_1$-isomorphic.
    \item These lens spaces are the types of a semi-free $\mathbb{Z}_n$-action  on a closed, oriented $4$-manifold $X$ with $m$ fixed points.
    \item There  exist  integers $k_1, ... , k_m$, each coprime to $n$, such that $\sum _{i=1}^m q_i k_i^2$ is divisible by $n$.
\end{enumerate}
{Moreover, if above conditions hold, then we can pick the manifold $X$ to be simply connected.}%\marginpar{Removed the simply-connected condition in (2) to get a full answer to Question 1.5 (2) in the cyclic case.}
\end{theorem}

The following theorem answers Question \ref{orbit type} (1).  %It is proved by applying Theorem \ref{main} to spherical 3-manifolds. 

\begin{theorem}\label{appl3}
For each spherical 3-manifold $Y$, there exists a closed, simply connected 4-manifold $X$ and an almost free $G$-action
with orbit type $Y$. Moreover, such an $X$ can be chosen such that the underlying space of $X/G$ is simply connected. 
\end{theorem}

\begin{remark}
The proof of Theorem \ref{appl3} can be adapted to any bounding spherical $n$-manifold for $n>3$.  Also note that there is no almost free action of $G$ on manifold $Y$ of dimension  $\leq 3$ such that the underlying space of $Y/G$ is simply connected \cite{Sc}.
\end{remark}

%\begin{theorem}\label{appl1}
%Two oriented lens spaces $L$ and $L'$ are $\pi_1$-isomorphic cobordism if and only if there is an degree one map $L\to L'$.
%\end{theorem}

%\begin{corollary}\label{appl2} For given two lens sapces $L, L'$,  there is an $\mathbb{Z}/n$ action on a simply connected closed 4-manifold with only $2$ fixed points  such that their types are  $L, L'$
 %if and only if they are homotopy equivalence.
%\end{corollary}

%We will prove  general versions of Theorem  \ref{appl1} and Corollary \ref{appl2} will be stated in  Theorem  \ref{appl1restated} and Corollary \ref{appl-2}.

\subsection{Complexity of 4-manifolds with given boundaries}
%{Finite index bounding and   achirality of 3-manifolds} 

Started from Hausmann-Weinberger \cite{HW},
some  3-manifold invariants are derived from related 4-manifolds, see \cite{SW1}  for more details. Given  Theorems \ref{main} and \ref{main3}, it is natural to consider the following new  invariant for bounding $n$-manifolds $Y$,  the minimal bounding index,  derived  from $(n+1)$-manifolds it bounds:
\[ 
\begin{split}
O_b(Y)=\min \{|\pi_1(W):\pi_{1}(Y)|\, \mid\, W \text{ is $\pi_1$-injectively bounded by }Y\}
\in \ZZ_+\cup\{\infty\}.    \end{split}\]
In particular $O_b(Y)$ is defined for each closed 3-manifold. 
%Intuitively, $O_{b}(Y)$ describes the minimal size  of a group that contains $\pi_{1}(Y)$ as the boundary subgroup. 
We say $Y$ is finite index bounding if $O_b(Y)< \infty$. %And we say $Y$ is index two bounding if $O_b(Y)=2$. 
%Clearly for each bounding $n$-manifold $Y$, $O_b(Y)$ can be defined,  
Clearly  $|\pi_1(Y)|<\infty$ implies  $O_b(Y)< \infty$ by Theorems 
\ref{main} and \ref{main3}.  
 
%We say that a 3-manifold $M$ {\it virtually} has a certain property if a finite cover of $M$ has this property. %which is based on Wise's machinery (\cite{Wise}). 
%The general philosophy is: a $3$-manifold has a lot of finite covers, so it should have a finite cover that satisfies any given reasonable property. 
 %see  a survey \cite{LS}.  We will address another two virtual properies of 3-manifolds: Virtually achiral and virtually bounds surface bundles over surfaces.

A closed orientable manifold is called {\em achiral},  if it admits an orientation reversing homeomorphism, and is called virtually achiral if it has an achiral finite cover.  The study of various virtual properties of $3$-manifolds became an active  topic on $3$-manifolds after Agol's solution (\cite{Ag}) of Thurston's virtual Haken conjecture \cite{Thu}.
%Virtually achirality of 3-manifolds is also discussed in a recent work \cite{TWWY}
%Both achirality  and virtual properties (such as virtually Haken and virtually special \cite{Ag}, \cite{PWi}) are  important topics in 3-manifolds and in geometric group theory.
The following results reveal some relations between  finite index bounding and virtual achirality and geometries of 3-manifolds:

\begin{theorem}\label{thm: O(Y)=2 or finite} Let $Y$ be a closed, orientable 3-manifold.

{\rm (1)}  If $Y$ is aspherical, then $O_{b}(Y)<\infty$ implies that  $Y$ is virtually achiral.
   
{\rm (2)}     If $Y$ admits an orientation reversing free involution, then $O_{b}(Y)=2$. The reverse is also true if $Y$ is hyperbolic. 

{\rm (3)}   Suppose $Y$ $\pi_1$-injectively bounds a compact orientable $4$-manifold $W$. Then for any integer $d>0$, $Y$ $\pi_1$-injectively bounds a compact orientable $4$-manifold $W_d$ such that $|\pi_1(W_d): \pi_1(Y)|=d|\pi_1(W): \pi_1(Y)|$.
\end{theorem}

%The next result claims that all odd primes are realized by len spaces

\begin{theorem}\label{O(Y)=3}
For each prime $p\ge 5$, $O_b(L(p,q))=\min\{d\ge 3|\, d|p-1\}.$ %Thus every odd prime is realized by $O_b(Y)$.
\end{theorem}

\begin{remark}\label{finite index} It is known that  {\rm (i)} each  $\mathcal G$ $3$-manifold is aspherical unless $ \mathcal  G$ is $S^2\times E^1$ or  $S^3$;
{\rm (ii)}  each  $\mathcal G$ $3$-manifold is not virtually achiral when  $\mathcal G$ is $Nil$ or $PSL(2,\mathbb{R})$. Moreover, \textcolor{red}{many} Sol and hyperbolic 3-manifolds are not virtually achiral \cite{TWWY}.
{\rm (iii)} there are $\mathcal G$ $3$-manifolds which admits orientation reversing free involution unless  $\mathcal G$ is either  $Nil$ or $\widetilde{PSL}(2,\mathbb{R})$, or $S^3$. 

{\rm (1)} By Theorem \ref{main3},  $O_b(Y)< \infty$ for each spherical $3$-manifold $Y$. By {\rm (i)}, {\rm (ii)} and Theorem \ref{thm: O(Y)=2 or finite},
$O_{b}(Y)=\infty$ for  each  $Nil$ or 
$\widetilde{PSL}(2,\mathbb{R})$ $3$-manifold $Y$.

{\rm (2)} If a closed orientable surface  $F$ $\pi_1$-injectively bounds
a compact orientable 3-manifold $Y$ with $|\pi_1(Y): \pi_1(F)|< \infty$, then $|\pi_1(Y): \pi_1(F)|=2$ \cite[Chap. 10]{He}.
 By {\rm (iii)} and Theorem \ref{thm: O(Y)=2 or finite} {\rm (2)}, for $\mathcal G\ne Nil, \widetilde{PSL}(2,\mathbb{R})$
and $S^3$, there exists  $\mathcal G$ 3-manifold $Y$ which $\pi_1$-injectively bounds a compact orientable 4-manifold $W_d$ with index $2d$ for any integer $d>0$. 

{\rm (3)}  By  \cite{Da}, $O_{b}(Y)=1$ if and only if $Y=S^3$ or a connected sum of $S^2\times S^1$. 
So any aspherical 3-manifold $Y$ has $O_b(Y)\ge 2$. Indeed  any aspherical n-manifold $Y$ has $O_b(Y)\ge 2$ \cite{SW2}. 

{\rm (4)} By {\rm (3)}, and by Theorem \ref{thm: O(Y)=2 or finite} and {\rm (iii)}, there are aspherical 3-manifolds $Y$ with $O_b(Y)=2$.
\end{remark}

% Inspired by Theorem \ref{O(Y)=3} and  Remark \ref{finite index},
%we wonder
%\begin{question} If each positive integer $n$ can be realized as the minimal bounding index?
%\end{question}

%More comprehensive  versions of Theorems \ref{thm: O(Y)=2 or finite},\ref{O(Y)=3} and Remark \ref{finite index} will appear in Section 4.

\subsection{Some explicit examples of 4-manifolds with $\pi_1$-injective boundaries} Except 3-manifolds described in Theorem \ref{thm: O(Y)=2 or finite} (2),
 it is usually hard to describe which and how  4-manifolds $W$ which are $\pi_1$-injectively bounded by given 3-manifolds $Y$.
%On the other hand, if we follow the developments  in 3-manifolds topology to discuss how $M$ virtually bounds $\pi_1$-injectively bounds a 4-manifold $W$, it become easier.
%{To see 4-manifolds bounded by 3-manifolds up to finite covers}
Surface bundles are important classes in both 3-manifolds and 4-manifolds. 
For 3-manifolds which are surface bundles, Proposition \ref{bundle} below 
  provides rather concrete description of
those bounded 4-manifolds $W$, which also has a flavor close to Question \ref{motivation}.

%Call a 3-manifold $M$ virtually  bounds a surfaces bundle over surface, if a finite cover of $M$ bounds a surfaces bundle over surface.
Let $\Sigma_g$ be the closed orientable surface of genus $g$. 
%It is clear that $M$ virtually  bounds a surfaces bundle over surface implies that $M$ is virtually surface fibered. 
%We observed that
%the inverse is also true. We prove that 
\begin{proposition}\label{bundle}  Suppose $Y$ is a $\Sigma_g$-bundle over $S^1$, $g\ge 3$. 
Then $Y$ bounds a surface bundle over a surface.
Moreover, the bounding is $\pi_1$-injective and $W$ has residually finite $\pi_1$.
\end{proposition}

%(1) is an observation based on the connections between the classifying space of the group oforiented homeomorphisms on $\Sigma_g$
%and $MCG_+(\Sigma_g)$, and (2) follow from (1),  Thurston's picture of 3-manifolds \cite{Thu},  
By Proposition \ref{bundle}
and  Agol and
Przytycki-Wise's virtual fibration results \cite{Ag}, \cite{PWi}, we have 
\begin{corollary}\label{AW}
Suppose $Y$ is a   closed orientable hyperbolic or mixed $3$-manifold. Then  a finite cover of $Y$ $\pi_{1}$-injectively bounds a surface bundle over surface.
\end{corollary}

%and some calculations  on $MCG_+(\Sigma_g)$ \cite{Mum}.

%we first show that any 3-manifold $Y$ has a finite cover $p:\widetilde{Y}\to Y$ that $\pi_{1}$-injectively bounds a 4-manifold $\widetilde{X}$. 

%For most cases, we set $\widetilde{Y}$ to be a surface bundle over a circle and set $\widetilde{X}$ to be a surface bundle over a surface. Existence of $\widetilde{Y}$ follows 

%And existence of $\widetilde{X}$ follows from the fact that the abelianization of a surface mapping class group is torsion.

By Theorem \ref{main},  for each spherical 3-manifold $Y$, we can define $\chi_b(Y)$ to be the minimum
$\chi(W)$  among all compact, orientable 4-manifolds $W$ with finite $\pi_1$ and $\pi_1$-injectively bounded by $Y$. %One can verify that
%$\chi_b(Y)= 1+ b_2(W) \ge 1$. Using the  achirality of some  lens space $Y$ 
We will explicitly construct 
some 4-manifold $W$ $\pi_1$-injectively bounded by $L(5,1)$ realizing both  
$O_b(L(5,1))=4$ and $\chi_b(L(5,1))=2$.

\section{Atiyah's generalization of Thom's Theorem and a surgery theorem}
\subsection{Results in dim $\ge 3$ for proving $\pi_1$-injective bounding results}

We use $H_n(X)$ to denote  $H_n(X, \ZZ)$ in the whole paper.

\begin{theorem}\label{Atiyah} Let $X$ be a CW-complex with $\tilde H_*(X)=0$. Let $M$ be 
a closed oriented $n$-manifold which is trivial in $\Omega_n$.
Then for any map $f: M\to X$, there is a compact oriented $(n+1)$-manifold $W$ such that  $\partial W=M$
and the map $ f: M\to X$ extends to a map $\tilde f: W\to X$.
\end{theorem}

\begin{proof}
The proof based on Atiyah's generalization of Thom's Theorem.

Recall that Atiyah defined the bordism homology group  $\{\MSO_k(X), k\ge 0\}$ and proved it is a generalized homology theory in \cite{Ati}. Let $R_k(X)=\MSO_k(X)$.
Let $c:X\to \{\text{point}\}$ be a constant map. It induces a map between Atiyah-Hirzebruch Spectral Sequence
$$c_*^2: E^2_{s, t}(X)\to E^2_{s,t}(\text{point})$$
where $E^2_{s, t}(X)=H_s(X, R_t(\text{point}))$ and $E^2_{s, t}(\text{point})=H_s(\text{point}, R_t(\text{point}))$. Since $\tilde H_*(X)=0$,
by universal coefficient theorem, $c_*^2$ is an isomorphism for all $s, t$.
Note that $$H_s(X, R_t(\text{point}))=H_s(\text{point}, R_t(\text{point}))=0$$ for all $s\ge 1, t\ge 0$. So 
$$ E^2_{s, t}(X)\to E^2_{s,t}(\text{point})=0$$
for all $s\ge 1, t\ge 0$. So the Atiyah-Hirzebruch Spectral Sequence collapses in $E^2$-page, so it collapses on $E^n$-page
for all $n\ge 2$. The Atiyah-Hirzebruch Spectral Sequence converges to bordism groups, it follows that the induced map $$c_* :    R_k(X)=\MSO_k(X)\to  R_k(\text{point})=\MSO_k(\text{point})$$
are isomorphism for all $k\ge 0$.

Reall each element in  $\MSO_k(X)$ is represented
by a map $f: M\to X$ where $M$ is a closed oriented $k$-manifold.
Then for any map $f:M\to X$ for a closed oriented $n$-manifold, consider the bordism class $[f: M\to X]$. Let 
$c_*[f: M\to X]=[c\circ f: M\to X]$ be the image in $\MSO_k(\text{point})$. Then it is represented by a map $c\circ f: M\to \{\text{point}\}$,
which is a constant map. Since $[M]=0\in \Omega_n$, there is a compact oriented $(n+1)$-manifold $W'$ such that  $\partial W'=M$.
Then $c\circ f$ extends to $W'$, that is,  there is a map $g:W'\to \text{point} $ with  $g|M=f$. 
It follows that $c_*[f: M\to X]$ is trivial in $\MSO_k(\text{point})$. Since $c_* :    \MSO_n(X)\to  \MSO_n(\text{point})$ is an isomorphism, we 
have that $ [f: M\to X]$ is trivial in $\MSO_n(X)$, that is there is a compact oriented $(n+1)$-manifold $W$ such that  $\partial W=M$
together with a map $\tilde f: W \to X$ which extends $f$.
\end{proof}

\begin{proposition}\label{surgery}
Suppose $\Gamma$ is a finitely presented group. Suppose $M$ is a closed oriented $n$-manifold, $n\ge 3$, and  $f: M\to K(\Gamma, 1)$ is a $\pi_1$-injective map. If $f$ extends to $\tilde f: W\to K(\Gamma, 1)$ for some  compact oriented $(n+1)$-manifold $W$
with $\partial W=M$, then we can choose $W$ so that $\tilde f$ is an $\pi_1$-isomorphism.
\end{proposition}

\begin{lemma}\label{finite generated} Let $\phi: G\to \Gamma$ be a surjection from a finitely generated group to a finitely presented group.
Then the kernel of $\phi$ is finitely normally generated.
\end{lemma}

\begin{proof}
Let $\phi: G\to \Gamma$ be a surjection from a finitely generated group to a finitely presented group.
Since $G$ is finitely generated, there is a surjection $\psi: F_n\to G$ from free group of rank $n$ for some $n$, therefore a surjection $\phi\circ \psi: F_n\to G \to \Gamma$.
Let $y_1, ... , y_n\in \Gamma$ be the images of the free generators $\{x_1, ... x_n\}$ of $F_n$ under $\phi\circ \psi$, then $y_1, ... , y_n$ is set of generators of $\Gamma$.
Since $\Gamma$ is finitely presented, and the property to be finitely presented is independent of the set of generators,
and we have a presentation
$$\Gamma=\langle y_1, ... , y_n\ \, | \, r_1(y_1, ..., y_n), ... , r_m(y_1, ..., y_n) \rangle,$$
which implies that the kernel of   $\phi\circ \psi$ is normally generated by 
$$\{r_1(x_1, ..., x_n), ... , r_m(x_1, ..., x_n) \}.$$
Then one can see directly  that the kernel of   $\phi$ is normally generated by 
$$\{\psi(r_1(x_1, ..., x_n)), ... , \psi(r_m(x_1, ..., x_n)) \}.$$
\end{proof}

\begin{proof}[Proof of Proposition \ref{surgery}]
Suppose $\tilde f: W\to K(\Gamma, 1)$ is an extension $f: M\to K(\Gamma, 1)$.
Let $k$ be the rank $\Gamma$.
Let $W_1=W\# (\#_k S^{n-1}\times S^1)$ be  the connected sums of $W$ and $k$ copies of $S^{n-1}\times S^1$.
Let $$\tilde f_1: W_1=W\# (\#_k S^{n-1}\times S^1) \to W\vee (\vee_k  S^1) \to K(\Gamma, 1)$$
be the composition of two maps: the first one pinch each $S^{n-1}\times S^1$ to $S^1$, and second one maps $W$ to $K(\Gamma, 1)$ via $\tilde f$, and maps those $k$
circles to the $k$ generators of $K(\Gamma, 1)$.
Clearly  $\tilde f_{1*}$ is  surjective
on $\pi_1$. 

Since $\tilde f_{1*}$ is a surjection between two finitely presented groups, by Lemma \ref{finite generated} the kernel of $\tilde f_{1*}$
is  normal generated by finitely many elements in $\pi_1(W_1)$. Let $c_1, ..., c_k$ be disjoint simple closed circles
in the interior of $W_1$ which represent the free homotopy classes of those generators.
Let $N(c_1), ... , N(c_k)$ be the disjoint regular neighborhood of $c_1, ... , c_k$ respectively. Then each 
$$N(c_i)\cong c_i\times D^n\cong S^1\times D^n.$$

Let $$W_2=W_1\setminus  (\cup_i c_i\times D^n)$$ and 
$$W_3=W_2\cup (\cup_i D^2_i\times S^{n-1}),$$
where each component $c_i\times S^{n-1}$ of $\partial W_2$ is identified with $\partial (D^2\times S^{n-1})=S^1\times S^{n-1}$ canonically.
Since $K(\Gamma, 1)$  has no homotopy groups of dimension $>1$, the restriction $\tilde f_1|: W_2 \to K(\Gamma, 1)$ extends to $\tilde f_3: W_3\to K(\Gamma, 1)$,
From Van Kampen theorem, it is easy to verify that $\tilde f_{3*}$ is an isomorphism on $\pi_1$.

Note during the surgery from $(\tilde f, W)$ to $(\tilde f_3, W_3)$, we do not touch $(f, M)$, we have a required extension $\tilde f_3: W_3\to K(\Gamma, 1)$.

\subsection{Results in dim $=3$ for further  applications} %In the text below, we often use $Y$ to denote closed oriented 3-manifolds

\begin{theorem}\label{pro: key proposition} Let $Y$ be a connected closed oriented $3$-manifold and let $\phi: \pi_{1}(Y)\to \Gamma$ be a group homomorphhism to a finitely presented group $\Gamma$. Let $f_{\phi}: Y\to K(\Gamma,1)$ be the map induced by $\phi$. Then the following two conditions are equivalent:

\begin{enumerate}
    \item There exists a smooth 4-manifold $X$ bounded by $Y$, and an isomorphism $\pi_{1}(X)\cong \Gamma$ under which $\phi$ is exactly the map induced by the inclusion $Y\to X$.
    \item The map 
       $f_{\phi,*}: H_{3}(Y;\mathbb{Z})\to H_{3}(K(\Gamma,1);\mathbb{Z})$
    is trivial. 
\end{enumerate}
\end{theorem}

%Our  proof of the main result, Theorem \ref{main}, relies the following Theorems \ref{Thom Atiyah} and \ref{virtual torsion} which will be proved in this section, as well as Propsositions \ref{fiber torsion}, \ref{H1=0} and \ref{surgery} which will be proved in next section.

\begin{theorem}\label{Thom Atiyah}%\cite{Ati}
 Suppose $X$ is a compact topological space,  $Y_1,...,Y_k$ are closed oriented 3-manifolds, and  $f_i:Y_i\rightarrow X$ are maps, $i=1,..., k$. If
$$\sum_{i=1}^{n}{(f_i)_*[Y_i]}=0,$$ Then there exists a $4$-manifold such that $\partial W=\cup_{i=1}^k Y_i$, and $f: W\rightarrow X$ such that $f|_{Y_i}=f_i$.
\end{theorem}

\begin{proof} We use the bordism homology groups $\MSO_k(X)$ \cite{Ati}. Consider the map $\psi_k: \MSO_k(X)\rightarrow H_k(X;\mathbb{Z})$ which sends $[Y,f]$ to $f_*[Y]$. It is known that $\MSO_*$ is a generalized homology theory.  (recall $\Omega_q=\MSO_q(\text{point})$) Thus there exists an Atiyah-Hirzebruch Spectral sequence whose $E^2$-page is $\{H_p(X,\Omega_q)\}$ and converges to $\{\MSO_k(X)\}$.  

Since $\Omega_q=0$ for $1\le q\le 3$, the spectral sequence collapses on $E_2$-page in the region $p+q\le 3$. Since $\Omega_0=\mathbb{Z}$, we have $E^2_{p,0}(X)=H_p(X, \mathbb{Z})$. So $E^\infty_{k,0}(X)=H_k(X,\mathbb{Z})$ for $k\le 3$.
Since $\Omega_q=0$ for $1\le q\le 3$, $E^\infty_{p,q}=0$ for $1\le q\le 3$. So we have $\MSO_k(X)\cong H_k(X,\mathbb{Z})$ when $k\le 3$.

We have the map $\psi_k: \MSO_k(X)\rightarrow H_k(X;\mathbb{Z})$ is an isomorphism when $k\le 3$.

Now consider the element $\xi=\sum_{i=1}^n [Y_i, f_i]\in \MSO_3(X)$. Then $\psi_3(\xi)=\sum_{i=1}^{n}{(f_i)_*[Y_i]}=0$. Since $\psi_3$ is an isomorphism, $\xi$ is trivial in $\MSO_3(X)$. It follows that there exists a $4$-manifold such that $\partial W=\cup_{i=1}^k Y_i$, and $f: W\rightarrow X$ such that $f|_{Y_i}=f_i$.
\end{proof}

\begin{proof}[Proof of Theorem  \ref{pro: key proposition}] (1) $\Longrightarrow$ (2):
Let $X\to K(\Gamma, 1)$ be the map which induces the identity on $\pi_1$. The the composition map $Y\to X\to K(\Gamma, 1)$
induces $\phi: \pi_1(Y)\to \Gamma$ on $\pi_1$, and the map  $f_{\phi,*}: H_{3}(Y;\mathbb{Z})\to H_{3}(K(\Gamma,1);\mathbb{Z})$
is trivial since the first map is trivial.

(2) $\Longrightarrow$ (1): Let $i=1$ and $Y=Y_1$,  by Theorem \ref{Thom Atiyah}, there exists a compact orientable 4-manifold $W$ such that 
$\partial W=Y$ and $\tilde f :W \to K(\Gamma, 1)$ such that $\tilde f|Y= f_{\phi}$. By Proposition \ref{surgery}, we can choose $W$
such that $\tilde f :W \to K(\Gamma, 1)$ is a $\pi_1$-isomorphism. So we have the commutative diagram

$$
\CD
Y     @>>>    W \\
@\vert @VV{\tilde f}V \\
Y   @>f_{\phi}>>K(\Gamma,1)
\endCD
$$
which induces commutative diagram on $\pi_1$
$$
\CD
\pi_1(Y) @>>> \pi_1(W) \\
@\vert @VV{\tilde f_*}V \\
\pi_1(Y) @>\phi>> \Gamma.
\endCD
$$

So under the isomorphism $\tilde f_*: \pi_1(W)\to \Gamma$, $\phi$ is exactly induced by the inclusion $Y\to W$.
\end{proof}

\section{Manifolds with (residually) finite $\pi_1$ bound manifolds with  (residually) finite $\pi_1$}

\subsection{A construction of finite mapping telescope $X_n$  keeping residual finiteness}

Let $X$ be a connected compact CW-complex, and choose a base point $x_0\in X$.
Let $$i_1, i_2: X\to X\times X$$ be given by
$i_1(x)=(x, x_0)$ and $i_2(x)=(x_0, x)$. 
Let $$\Delta: X\to X\times X$$ be the diagonal embedding given by
$\Delta(x)=(x, x)$. 

\begin{figure}[h]
\includegraphics[width=240pt, height=140pt]{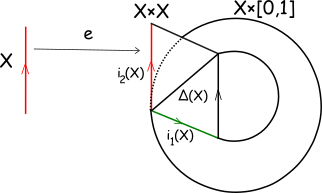}
\caption{Sketch picture for $\sigma(X)$}\label{fig:figure1}
\end{figure}

%Then define $\sigma(X)$ to be the colimit of $X\to X\times X$ via the embeddings $i_1$ and $\Delta$. Precisely 
Let $\sigma(X)$ be the quotient space
$$\sigma(X)=\frac{X\times X\coprod X\times [0,1]}\sim,$$
where  $\Delta(X)$ 
is identified with $X\times \{0\}$ via $(x, x)\sim (x, 0 )$, and $X\times x_0$ is identified with $X\times \{1\}$ via $(x, x_0)\sim (x, 1 )$. See Figure 1 for sketch picture of $\sigma(X)$. Let 
$$q: X\times X\to \sigma(X) $$ be the quotient map.

Since $X$ is compact, $\sigma(X)=X\times X/\sim$ is also compact.
Moreover since $X$ is a CW-complex, so is $\sigma(X)$.
Consider the composition
$$e=q\circ i_2: X\to X\times X\to \sigma(X), \qquad (*)$$
which is an embedding from $X$ to $\sigma(X)$. We will repeat this construction several times in our argument.

For each group $G$ with unit $1$, if we define $\sigma(G)$ to be an HNN extension of $G\times G$ by $t$:
$$\sigma(G)=\langle G\times G, t \, | \, t(g, 1) t^{-1}=(g,g), \, \text{for any}\, g\in G\rangle,$$
There is also a homomorphism of groups $$e=\beta\circ i_2: G\to G\times G\to \sigma(G),$$
where $i_2(g)=(1,g)$ and $\beta : G\times G\to \sigma(G)$ is the canonical inclusion \cite{ScW}. By Van Kampen theorem, one can verify the following result.

\begin{lemma}\label{G} The fundamental group of $\sigma(X)$ is given by
$$\pi_1(\sigma(X))=\sigma(\pi_1(X)).$$
Moreover the induced map of the embedding $e: X\to \sigma (X)$ is exactly the homomorphism 
$$i: \pi_1(X)\to \sigma(\pi_1(X))=\pi_1(\sigma(X))$$
defined above.
\end{lemma}

%Let $\sigma(X)$ be the quotient space
%$$\sigma(X)=X\times X\prod /\sim,$$
%where $i_1(x)=(x, x_0)\sim (x,x )=\Delta(x)$ for each $x\in X$. Let 
%$$q: X\times X\to \sigma(X) $$ be the quotient map. 
%Then since $X$ is compact, $X\times X$ is compact. Therefore $\sigma(X)=X\times X/\sim$ is compact.
%Moreover since $X$ is a CW-complex, so is $\sigma(X)$.

For each connected compact CW-complex  $X$, we define a sequence of spaces and embeddings as below:
Let $X=X_0$ and let $X_{n}=\sigma^{n}(X)$. Then$X_n=\sigma(X_{n-1}).$ 
Then we have the embedding 
$$e_n: X_n\to X_{n+1}$$ given by (*).
Now
the mapping telescope $X_\infty$
of the embedding sequence 
\begin{equation}\label{eq: sequence Xn}
X_0\to X_1\to X_2\to ... \to X_{n-1} \to X_n\to ...    
\end{equation}
is defined as 
$$X_\infty=\bigsqcup X_n\times [0,1]/\sim,$$
  where $(x_n, 1)\sim (x_{n+1}, 0)$ if  $e_n(x_n)=x_{n+1}$.

%Then given $n$, we have an embedding
%$$\tau_n=e_{n-1}\circ .... \circ e_1\circ e_0 :X\to X_1\to X_2\to ...  \to X_{n-1} \to X_n. \qquad (**)$$

\begin{proposition}\label{tech3} Suppose $G$ is finitely generated group. Then 

{\rm (1)} $e: G\to \sigma(G)$ is injective.

{\rm (2)} $\sigma(G)$ is residually  finite if $G$ is.
\end{proposition}

%A group $G$ is said to be residually finite if for each $1\ne g\in G$, there is a finite goup $H$ and a homomorphism $\phi: G\to H$
%such that $\phi(g)\ne 1$. 

The proof of Proposition \ref{tech3} (2) need more explicit description of HNN extension and  some results.
Given a group $\Gamma$, subgroups $C_0, C_1$, and an isomorphism $\phi: C_0\to C_1$,
  we have the so called HNN extension $\Gamma$ by identifying $C_0$ and $C_1$ vis $\phi$,
  denoted as $\HNN(\Gamma, C_0, C_1, \phi)$ \cite{ScW}, \cite[Chap. 15]{He}. Then we have 
  
  $$\sigma(G)=\HNN(G\times G, C_0, C_1, \phi),  $$ 
  where  $C_0=\{(g, g)| g \in G\}$, $C_1=\{(g, 1)| g\in G\}$, and 
  $\phi : C_0\to C_1$ is given by $\phi((g, g))=(g, 1)$.
  
 %The proof of that $\sigma(G)$ is residually finite is based on the following facts:

 \begin{proposition} \label{4} \cite[15.20. Lemma]{He}
 Let  $H=\HNN(\Gamma, C_0, C_1, \phi)$ with finitely generated $\Gamma$, $C_0$ and $C_1$.
 Suppose there is a sequence $\{N_i\}$ of normal subgroups of finite index in $\Gamma$ satisfying

 {\rm (i)} $\cap N_i=1$,
 
 {\rm (ii)} $\cap N_iC_0=C_0$, $\cap N_iC_1=C_1$, and
 
 {\rm (iii)} $\phi(N_i\cap C_0)=N_i\cap C_1$ for all $i$.
 
 Then $H$ is residually finite.
 \end{proposition}
 
 \begin{lemma}\label{5} \cite[15.16. Lemma]{He}
For a finitely generated group $G$, $G$ is residually finite if and only if the intersection 
of all its finite index subgroups is trivial.
\end{lemma}

\begin{proof}[Proof of Proposition \ref{tech3}]
(1)  Recall $e=\beta\circ i_2: G\to G\times G\to \sigma(G),$
where $i_2(g)=(1,g)$  clearly is injective, and  the canonical map $\beta : G\times G\to \sigma(G)$ is also injective \cite[Theorem 1.7]{ScW}. So $e$ is injective.

(2) Since $G$ is finitely generated, all $G\times G$, $C_0=\{(g, g)| g \in G\}$ and  $C_1=\{(g, e)| g\in G\}$
are finitely generated. 

Since $G$ is residually finite, there is a sequence $\{K_i\}$ of normal subgroups of finite index in $G$ satisfying
$ \cap K_i=1$ by Lemma \ref{5}. Let $N_i=K_i\times K_i$,  it is easy to see that the  $\{N_i\}$ is a sequence of normal subgroups of finite index in $\Gamma$ satisfying $ \cap N_i=1$, that is, the condition (i) in Proposition \ref{4} is satisfied.

Next we verify the condition (ii) in Proposition \ref{4} is satisfied. We just verify that $\cap N_iC_0=C_0$. Clearly $C_0\subset 
\cap N_iC_0$. On the other hand, we have
$$N_iC_0=(K_i\times K_i)C_0=\{(k_i, k'_i)(g,g)|k_i, k'_i\in K_i, g\in G \}$$
$$=\{(k_ig, k'_ig)|k_i, k'_i\in K_i, g\in G\}=\{(g_1, g_2)| g_1g_2^{-1}\in K_i\}.$$

Suppose $z\notin B_0$, then $z=(g_1, g_2)$ such that $g_1g_2^{-1}\ne 1$. Since $ \cap N_i=1$, $g_1g_2^{-1}\notin K_i$ for some $i$,
that is $z\notin N_iB_0$ for some $i$. We finish the verification of (ii).

Finally we verify the condition (ii) in Proposition \ref{4} is satisfied. Note
$$ N_i\cap C_0= \{(g, g)| g\in K_i\}, \, \, N_i\cap C_1= \{(g, 1)| g\in K_i\}. $$
Then clearly $z\in  N_i\cap C_0$ if and only if $\phi(z)\in  N_i\cap C_1$. We finish the verification of (iii).

Therefore $\sigma(G)$ is residually finite.
\end{proof}

%\begin{proposition}\label{vr3}
%The embedding $\tau_n: X\to X_n$ is $\pi_1$-injective.
%\end{proposition}

%\begin{proposition}\label{rf3}
%If $\pi_1(X)$ is residually finite, so is $\pi_1(X_n)$.
%\end{proposition}

%Proposition \ref{zero1},  Proposition \ref{vr3}  and Proposition \ref{rf3}  will be proved in Section 5,  Section 6 and Section 7 respectively.

\subsection{The infinite mapping telescope $X_\infty$ with trivial homology}

\begin{proposition}\label{zero1}
$H_*(X_\infty)=H_*(\text{point})$
\end{proposition}

Consider the sequence (\ref{eq: sequence Xn}). For $n>m$, we define the map  
$$\tau_{m, n}= e_{n-1}\circ ... \circ e_m: X_m\to...  \to X_{n-1}\to X_n.$$
Then we have the following property of $\tau_{m,n}$ on homology groups.

\begin{lemma} The following are equivalent:

{\rm (1)} For any integers $d>0$, and $N>0$, there exists an $n>N$ such that $\tau_{N, n}: X_N\to X_n$ induces trivial maps on $H_i$ for $1\le i \le d$.

{\rm (2)} $\tilde H_i(X_\infty)=0$.
\end{lemma}
 
 \begin{proof}
 Suppose (1) holds. For any $k$-cycle $c\in X_\infty$, $c\subset X_N$ for some $N$. Then for some $n>N$, $c=\partial D$ for some
 $D\subset X_n \subset X_\infty$. Hence $c$ is zero in $H_*(X_\infty)$, i.e. $\tilde H_i(X_\infty)=0$.
 
 Suppose (2) holds. For each $i\in \{1, ... ,  d\}$, fix a finite generating set of $H_i(X_N)$. For  any element $c$ in this basis, since $\tilde H_i(X_\infty)=0$, $c=\partial D$ for some finite chain $D\subset X_\infty$.
 Since $D$ is compact, $D\subset  X_{n_c}$. 
 Since there are only finitely many elements in this set, there exists an $n_i>N$ such that each element in this set bounds in $X_{n_i}$. Then the image of $H_i(X_N)$ vanishes in $H_i(X_{n_i})$.
 Let $n=\text{max}\{n_i, i=1, ..., d\}$, we have
 that $\tau_{N, n}: X_N\to X_n$ is trivial in $H_i$ for $1\le i \le d$.
 \end{proof}
 
 So to prove (2), we need only to prove (1), and to prove (1), we need only to prove the following
 
 \begin{proposition}\label{1}
$X\to \sigma^{3^{n-1}}(X)=X_{3^{n-1}}$ induces trivial maps on $H_i$ for $1\le i\le n$.
\end{proposition}

We will prove Proposition \ref{1} by induction based on the following

 \begin{proposition}\label{2}
Suppose we have a composition  
$$A_1\overset {f_1}  \longrightarrow A_2\overset {f_2}  \longrightarrow  A_3 \to \sigma (A_3).$$
If $f_1$ and $f_2$ induce  trivial maps  on $H_i$ for  $1\le i\le n-1$, then the composition $A_1\to \sigma(A_3)$ induces  trivial maps in $H_i$ for  $1\le i\le n$.
\end{proposition}

To start the induction, we need 

\begin{lemma}\label{3}
$X\to \sigma(X)$ induces trivial map on $H_1$.
\end{lemma}

\begin{proof} Recall  $$i_1: X\to X\times X, \, i_2:  X\to X\times X, \, \Delta : X\to X\times X$$
 be the embedding of $X$ to the first factor, the second factor,
and diagonal map respectively, and  
$$q: X\times X \to \sigma(X)$$
be the quotient map.
Since in  construction of $\sigma(X)$, the first factor $X$ and the diagonal are identified canonically,  we have
$$q\circ\Delta=q\circ i_1  $$
and the embedding $e: X\to \sigma (X)$ is given by 
$$e=q\circ i_2.$$
Applying Kunneth formular \cite[Theorem 3B.6.]{Ha1} to $H_1(X\times X)$, since $\textrm{Tor} (H_0(X), H_0(X))=\textrm{Tor}(\mathbb{Z},\mathbb{Z})=0$,
we have $$H_1(X\times X)=H_1(X)\otimes \mathbb{Z} \oplus \mathbb{Z}\otimes  H_1(X\times X).$$

Then one can derived that

$$i_{1*}+i_{2*}=\Delta_*.$$
So we have 
$$e_*=q_*\circ i_{2*}= q_*\circ (\Delta_*-i_{1*})= q_*\circ \Delta_*- q_*\circ i_{1*}=(q\circ \Delta)_*-( q\circ i_1)_*=0,$$
that is,  the embedding induces trivial map on $H_1$.
\end{proof}

\begin{proof}[Proof of Proposition \ref{1}] By Lemma \ref{3}, Proposition \ref{1} hold for $k=1$.

Suppose Proposition \ref{1} hold for $k=n-1$.
Consider the embedding sequence
$$X\to \sigma^{3^{n-1}}(X)\to  \sigma^{ 3^{n-1}}(\sigma^{3^{n-1}}(X))=\sigma^{2\times 3^{n-1}}(X)$$
$$\to \sigma(\sigma^{2\times 3^{n-1}}(X))= \sigma^{2\times 3^{n-1}+1}(X)
\to \sigma^{3^{n}}(X).$$
By the induction hypothesis on $n-1$, the first two maps induce trivial maps on $H_i$ for $1\le i\le n$.
By Proposition \ref{2}, the embedding 
$$X\to \sigma(\sigma^{2\times 3^{n-1}}(X))= \sigma^{2\times 3^{n-1}+1}(X)$$ induces trivial maps on $H_i$ for $1\le i\le n+1$, therefore 
the embedding 
$$X\to \sigma^{3^{n}}(X)$$ induces trivial maps on $H_i$ for $1\le i\le n$.
\end{proof}

\begin{proof}[Proof of Proposition \ref{2}]
 We will prove that for the sequence
$$A_1\overset {f_1}  \longrightarrow A_2\overset {f_2}  \longrightarrow  A_3 \overset {e}  \longrightarrow\sigma (A_3).$$
if$f_1$ and $f_2$ induce  trivial maps  on $H_i$ for  $1\le i\le n-1$,  then composition $A_1\to \sigma(A_3)$ are trivial in $H_i$ for  $1\le i\le n$.

\end{proof}

We start from the following commutative diagram

$$\begin{CD}
			A_1 @>  f_1 >> A_2 @>  f_2 >> A_3   \\
		         @VV \Delta_1 V   @VV \Delta_2 V   @ VV \Delta_3V\\
			A_1\times A_1@> f_1\times f_1 >> A_2\times A_2@> f_2\times f_2 >>A_3\times A_3.
		\end{CD} \qquad (1)  $$
		
		Then we have the  following commutative diagram in $H_n$
		
		$$\begin{CD}
			H_n(A_1) @>  f_1 >> H_n(A_2) @>  f_2 >> H_n(A_3)   \\
		         @VV \Delta_1 V   @VV \Delta_2 V   @ VV \Delta_3V\\
			H_n(A_1\times A_1)@> f_1\times f_1 >> H_n(A_2\times A_2)@> f_2\times f_2 >>H_n(A_3\times A_3).
		\end{CD} \qquad (2)  $$

		Apply Kunneth formula \cite[Theorem 3B.6.]{Ha1} to the second low of (1), we have the following commutative diagram
		
		$$\begin{CD}
			0 @>  >>\overset\bigoplus {_{k+l=n}} H_k(A_1)\otimes H_l(A_1) @>  j_1 >> H_n(A_1\times A_1) @>  p_1 >> \overset\bigoplus {_{k+l=n-1}} \text{Tor}( H_k(A_1), H_l(A_1))   \\
			@ VV V@VV  f_1\otimes f_1 V   @VV f_1\times f_1 V   @VV  V   \\
			0 @>  >>\overset\bigoplus {_{k+l=n}} H_k(A_2)\otimes H_l(A_2)@>  j_2 >> H_n(A_2\times A_2) @>  p_2 >> \overset\bigoplus {_{k+l=n-1}} \text{Tor}( H_k(A_2), H_l(A_2))\\
			@ VV V@VV f_2\otimes f_2 V  @VV  f_2\times f_2  V   @VV  V   \\
			0 @>  >>\overset\bigoplus {_{k+l=n}} H_k(A_3)\otimes H_l(A_3)@>  j_3 >> H_n(A_3\times A_3) @>  p_3 >>\overset\bigoplus {_{k+l=n-1}} \text{Tor}( H_k(A_3), H_l(A_1)).
		\end{CD}   \qquad(3)$$ 
		
		For each $\alpha \in H_n(A_1)$, we are going to prove $e \circ f_2 \circ f_1(\alpha)=0.$
		 
		 Set $\alpha_2=f_1(\alpha)$ and $\alpha_3=f_1(\alpha_2)$. 
		 
		 Now we explain the roles of $f_1$ and $f_2$ in Proposition \ref{2}: $f_1$ is to ensure $\Delta_2(f_1(\alpha))$ projects to $0\in \oplus_{i+j=n-1}{\Tor}(H_i(A_2),H_j(A_2))$, therefore  it is an image of an element  $\tilde \alpha_2\in  \oplus_{i+j=n}H_i(A_2)\otimes H_j(A_2)$; $f_2$ is to ensure the image of $\tilde \alpha_2$ in $\oplus_{i+j=n}H_i(A_3)\otimes H_j(A_3)$ has only  component
		 with $i=0$ or $j=0$.

		By conditions posed on $f_1$, the right-up vertical homomorphism  in the above diagram is trivial.
		Then by using the commutativity of the right-up square of (3), we have $p_2\circ (f_1\times f_2)=0$.
		So we have 
		$$0=p_2\circ (f_1\times f_2) \circ \Delta_1(\alpha)=p_2 \circ \Delta_2  \circ f_1(\alpha)=p_2 \circ \Delta_2(\alpha_1)$$
where the second $"="$ comes from the 	the commutativity of the left square of (2). So $\Delta_2(\alpha_1)\in ker (p_2)$.
By the exactness of second low of (3), we have $j_2(\tilde \alpha_2)=\Delta_2(\alpha_1)$ for some 
\[\tilde \alpha_2\in \overset\bigoplus {_{k+l=n}} H_k(A_2)\otimes H_l(A_2).\]
Let $\tilde \alpha_3=f_2\otimes f_2(\tilde \alpha_2)$. By conditions posed on $f_2$, we have
$$\tilde \alpha_3=\alpha_{n,0}+ \alpha_{0,n},$$
where $\alpha_{n,0}\in    H_n(A_2)\otimes H_0(A_2),\,    \alpha_{0,n}\in  H_0(A_2)\otimes H_n(A_2)$.

From the definition (or construction) of $j_3$, we have
$$j_3(\tilde \alpha_3)=j_3( \alpha_{n, 0})+ j_3( \alpha_{0,n})=i_1(\beta_1)+i_2(\beta_2),$$
where $\beta_1, \beta_2\in H_n (A_3)$. Then we have

$$\Delta_3(\alpha_3)=j_3(\tilde \alpha_3)=i_1(\beta_1)+i_2(\beta_2),$$
where the first $"="$ follows from the commutativity of both the right square of (2) and middle-down square of (3).
Let $$p_1: A_3\times A_3\to A_3$$ be the projection  to its the first factor, we have
$$p_1\circ i_1=id_{A_3}, \, p_1\circ \Delta=id_{A_3}, \, p_1\circ i_2=0.$$

So 
$$\alpha_3=(p_1\circ \Delta)(\alpha_3)=p_{1}\circ (i_1(\beta_1)+i_2(\beta_2))=p_{1}\circ i_1(\beta_1)+p_{1}\circ i_2(\beta_2)=\beta_1.$$
Similar arguments show that $\alpha_3=\beta_2$. So we have

$$\Delta_3(\alpha_3)=i_1(\alpha_3)+i_2(\alpha_3).$$

Now consider the quotient map $q_3: A_3 \times A_3\to \sigma(A_3)$. As we see in the proof of Lemma \ref{3},
$q_3\circ  \Delta_3(\alpha_3)=q_3\circ  i_1(\alpha_3)$, so we have 
$$0=q_3\circ  i_2(\alpha_3)=q_3\circ  i_2 \circ  f_2\circ f_1(\alpha)=e \circ  f_2\circ f_1(\alpha).$$
This finishes the proof.
\end{proof}

\subsection{Manifolds with  finite $\pi_1$ bound manifolds with  finite $\pi_1$} In this section, we prove Theorem \ref{main}. We start with a 
algebraic lemma.

\begin{lemma}\label{Reduction Thm}
Suppose $G$ is a finitely-generated residually-finite group and $H$ is a finite group and $\phi:H\to G$ is an injective homomorphism.
Then there exists a finite group $G_1$ and a homomorphism $\psi:G\to G_1$ such that the composite map
$$\psi\circ\phi:H\to G\to G_1$$
is injective.
\end{lemma}
\begin{proof}
Note that $\phi(H)\subset G$ is a finite subgroup. Since $G$ is residually-finite, for any $h\in H$, $h\ne e$, there exists a finite-index normal subgroup $N(h)\subset G$ such that $h\notin N(h)$.
Write $H=\{h_1=e,h_2,...,h_m\}$ where $m=|H|$. Then for any $2\le i\le m$, there exists a finite-index normal subgroup $N_i\subset G$ such that $h_i\notin N_i$.
Let $N=\cap_{i=2}^mN_i$. Then $N\subset G$ is a finite-index normal subgroup.
Note that for any $2\le i\le m$, we have $h_i\notin N_i$. So $h_i\notin N$. Therefore $H\cap N=\{e\}$.
 Let $G_1=G/N$ and $\psi:G\to G_1$ be the quotient map.
Then
$$ker(\psi)\cap \phi(H)=N\cap \phi(H)=\{e\}.$$
Therefore we get an  injective composite map
$$\psi\circ\phi:H\to G\to G_1.$$
\end{proof}

Now we restate Theorem \ref{main} as Theorem \ref{main2}.  Theorem \ref{main2} (1) is known  \cite{Hau}, we reprove it in our route,
then use it to prove  Theorem \ref{main2} (2) and (3), that is our Theorem \ref{main}.
\begin{theorem}\label{main2}
Suppose $M$ is a closed oriented bounding $n$-manifold. Then 
\begin{enumerate}
    \item $M$ $\pi_1$-injectively bounds a compact oriented $(n+1)$-manifold.
    \item $M$ $\pi_1$-injectively bounds a compact oriented $(n+1)$-manifold
with residually finite $\pi_1$ if $\pi_1(M)$ is residually finite.
\item $M$ $\pi_1$-injectively bounds a compact oriented $(n+1)$-manifold with finite $\pi_1$ if $\pi_1(M)$ is finite. 
\end{enumerate}
 \end{theorem}

\begin{proof}[Proof of Theorem \ref{main2}]
(1) Let $X_0=M$ be a closed oriented $n$-manifold.  Then we have  the  sequence of embeddings and its mapping telescope
$$\tau: X=X_0\to X_1\to X_2\to ...  \to X_{n-1} \to X_n\to ...\to  X_\infty.$$

By Proposition \ref{zero1}, we have $H_*(X_\infty)=H_*(\text{point})$. Since $M=0\in \Omega_n$,
by Theorem \ref{Atiyah}, the map $\tau: M\to X_\infty$ extends to a map $\tilde \tau: W\to X_\infty$ for a compact $(n+1)$-manifold
$W$ such that $\partial W=M$, more precisely
$$\tilde \tau \circ i =\tau : M\to X_\infty,$$
where $i: M\to W$ is the inclusion.
Since $W$ is compact, $\tilde \tau(W)\subset X_n$ for some $n$. 
By Proposition \ref{surgery}, we may assume the inclusion map $\tilde\tau_n: W\to X_n$ is $\pi_1$-isomorphic.
Then we have 
$$\tilde \tau_n \circ i =\tau_n : M\to X_n.$$

By Proposition \ref{tech3} (1), $\tau_n: M\to X_n$ is $\pi_1$-injective. 
Since $\tilde \tau_n$ is $\pi_1$-isomorphism, it concludes that $i: M\to W$ is $\pi_1$-injective.

(2) If $\pi_1(M)$ is residually-finite, then by Proposition \ref{tech3} (2), $\pi_1(X_n)$ is residually-finite. Since $\pi_1(W)\cong\pi_1(X_n)$, we get that $\pi_1(W)$ is residually-finite. 

(3) Suppose $\pi_1(M)$ is finite. Then it is residually-finite. By (2), $M$ $\pi_1$-injectively bounds 
a compact oriented $(n+1)$-manifold $W_0$ with residually-finite $\pi_1$.

Let $i_0:M\to W_0$ be the inclusion map. Then we have an injective map
$$\phi=(i_0)_*:H=\pi_1(M)\to\pi_1(W_0)=G.$$
%Then $H$ is a finite group and $\phi:H\to G$ is injective.
Now apply Lemma \ref{Reduction Thm}, there is a finite group $G_1$ and a homomorphism $\psi:G\to G_1$ such that the composite map
$$\psi\circ\phi:H\to G\to G_1$$
is injective.
There exists a map
$$F:W_0\to K(G_1,1)$$
such that $F_*=\psi:\pi_1(W_0)\to G_1.$
Let $$f=F|_M=F\circ i_0:M\to W_0 \to K(G_1,1).$$
Clearly $f$ extends to $W_0$ and $f_*=\psi\circ\phi:\pi_1(M)\to G_1$ is injective.
By Proposition \ref{surgery}, there exists another compact oriented $(n+1)$-manifold $W$ with $\partial W=M$ such that $f$ can be extended to $F':W\to K(G_1,1)$
such that the induced map
$$F'_*:\pi_1(W)\to K(G_1,1)$$ is an isomorphism. 
Let $i:M=\partial W\to W$ be the inclusion map. Since $$f=F'|_M=F\circ i:M\to W\to K(G_1,1),$$ and $f$ is $\pi_1$-injective, we get the inclusion map $i:M\to W$ is $\pi_1$-injective. Note that $\pi_1(W)\cong G_1$ is finite. 
\end{proof}

\section{Finite group actions on 4-manifolds and $\pi_1$-isomorphic cobordism lens spaces}

%\begin{theorem}\label{main}
%Each closed orientable  3-manifold $Y$ $\pi_1$-injectively bounds a compact orientable 4-manifold $X$. Moreover we can require $|\pi_1(X)|\le \infty$
%if $|\pi_1(Y)|\le \infty$.
%\end{theorem}

In this section we will prove Theorem \ref{appl-2}, Theorem \ref{appl1} and Theorem \ref{appl3}. 
Let $S^3$ be the unit sphere of $\mathbb{C}^2$.
Define a cyclic group action  $\tau_{p,q}: \mathbb{C}^2\to \mathbb{C}^2$ by
$\tau_{p,q}: (z_1,z_2)\mapsto(e^{\frac{2\pi i}{p}}z_1,e^{\frac{2q\pi i}{p}}z_2).$
Then for each pair of coprime integers $(p, q)$, $p>0$, we have  $L(p, q)= S^3/\tau_{p,q}.$
Now $S^3$ has the induced orientation from the unit $4$-ball $B^4\subset\mathbb{C}^2$ and $L(p,q)$ has the induced  orientation from the covering $S^3\to L(p,q)$.

\subsection{$\pi_1$-isomorphic cobordisms of lens spaces and semi-free $Z_n$-actions} 

The following theorem is a slight refinement of Theorem \ref{appl1}.

\begin{theorem}\label{appl1restated} Let  $L(n, q_1), ... , L(n, q_m)$  be $m$ oriented lens spaces.  Then the following conditions are equivalent:

\begin{enumerate}
    \item There is a compact, oriented, connected  4-manifold $W$ such that  $\partial W= \bigcup _{i=1}^m L(n, q_i)$ and each inclusion $L(n, q_i)\to W$ is $\pi_1$-isomorphic.
    \item These lens spaces are exactly the types of a semi-free  $Z_n$ action  on a closed oriented connected 4-manifold $X$ with $m$ fixed points.
    \item There  exist integers $k_1, ... , k_m$, each coprime to $n$, such that $ \sum _{i=1}^m q_i k_i^2$ is divisible by $n$.
    \item There is a $\pi_1$-isomorphic map $g_i: L(n, q_i) \to L(n,1)$ for each $i$ such that $\sum _{i=1}^m \degr(g_i)=0$.
\end{enumerate}
Moreover, we can pick the manifold $X$ in (2) to be simply-connected.
\end{theorem}

\begin{proof} 
(1) $\Longrightarrow$ (4):  Let $L=L(n,1)$ and $L_i=L(n, q_i)$.
Suppose first there is an oriented compact 4-manifold $W$ 
such that $\partial W= \bigcup _{i=1}^m L_i$ and each inclusion $L_i\to W$ is $\pi_1$-isomorphic.

Since $\pi_2(L)=0$, we can build a $K(\pi_1(L), 1)$ space $K$ by attaching cells of dimension $>3$ to $L$.
So there is an embedding $e: L\to  K$ as the 3-skeleton.
Then $H_3(K)=Z_n$ and $e_*=\operatorname{id}$ on $\pi_1$. Moreover, $e_{*}[L]\in H_3(K)$ is a primitive element.

Let $e_i: L_i\to W$  be the inclusions for $i=1..., m$. Then 
$$\sum _{i=1}^m e_{i*}[L_i]=0\in H_{3}(W).$$ 

Since $\pi_1(W)=\mathbb{Z}/n$, there is a $\pi_1$-isomorphic map $f: W \to K$.
%It follows that $f_*(\alpha) =f_*\circ j_{0\#}([L])=i_*[L]$ is primitive in $H_3(K(\pi_1(L), 1)$.
Now consider the map $f\circ e_i: L_i\to W \to K$. By cellular approximation theorem, there is map $f_i: L_i\to L$
such that $f\circ e_i$ is homotopic to $e\circ f_{i}$. Since both $f$ and $e_i$ are $\pi_1$-isomorphic, so is $e_i$. So we have 
\begin{equation}\label{eq: degree sum}
0=f_*(\sum_{i=1}^m e_{i*}[L_i])= \sum _{i=1}^m e_* f_{i*}([L_i])
=\sum _{i=1}^m e_*(\degr (f_i)[L])=(\sum _{i=1}^m \degr (f_i))\cdot e_*([L])    
\end{equation}
{Since $e_{*}[L]$ is primitive in $H_3(K)=\mathbb{Z}/n$, (\ref{eq: degree sum}) implies that $\sum _{i=1}^m\degr (f_i)=kn$ for some integer $k$.
 Let $g_{1}:L_{1}\to L_{1}$ be the composition
 \[
 L_1\cong L_1\#S^{3}\xrightarrow{q} L_{1}\vee S^{3}\xrightarrow{ f_{1}\vee p_{-kn}} L_{1}.  
 \]
 Here $q$ is the map that pinches the 2-sphere in the connected sum to a point, and 
$p_{-kn}: S^3\to L_1$ is a map of degree $-kn$.  So we have  $\degr(g_1)=\degr (f_1)-kn$. 
Now let $g_i=f_i$ for $i=2,...,n$. we have $\sum _{i=1}^m\degr (g_i)=0$.}

(4) $\Longrightarrow$ (1): Suppose that there are  $\pi_1$-isomorphic maps
\[g_i: L_i \to L \text{ for }1\leq i\leq m\] 
such that $\sum _{i=1}^m \degr(g_i)=0.$
%Then there is a map $$\bigcup_{i=1}^m g_i: \bigcup_{i=1}^m L_i\to L$$
Then we have  $\sum_{i=1}^m g_i([L_i])=0$ in $H_3(L)$. Then by Theorem \ref{Thom Atiyah},
there is compact oriented 4-manifold $W$ with $\partial W=\bigcup_{i=1}^m L_i$ and a map $\title f: W\to L$ which extends the map  
$$\bigcup_{i=1}^m g_i: \bigcup_{i=1}^m L_i\to L.$$  Moreover we can require that $\widetilde f: W \to L$ is a $\pi_1$-isomorphism by Proposition
\ref{surgery}. Then the inclusion of each  $L_i\to W$ is $\pi_1$-isomorphic .

(4) $\Longrightarrow$ (3): Suppose there exists a $\pi_1$-isomorphic map \[g_i: L(n,q_i)\rightarrow L(n,1)\] for each $i$ such that $\sum_{i=1}^m{\degr(g_i)}=0$. By Lemma \ref{degree} below, we have  \[\degr(g_i)=q_i k_i^2+n x_i\] for some $k_i$ coprime to $n$ and some $x_i\in \mathbb{Z}$. Hence we have
$$\sum_{i=1}^m q_i k_i^2 =-n\sum_{i=1}^m x_i.$$

(3) $\Longrightarrow$ (4): Suppose there exist {$k_1,...,k_m$}, each coprime to $n$, such that $\sum_{i=1}^m q_i k_i^2=xn$ for some integer $x$.
Let $x_1=-x, x_2=...=x_m=0$. Then there exists a $\pi_1$-isomorphic map $g_i: L(n,q_i)\rightarrow L(n,1)$ with degree $q_i k_i^2+nx_i$ by Lemma \ref{degree}.
Then
$$\sum_{i=1}^m\degr(g_i)=\sum_{i=1}^m(q_i k_i^2+n x_i)=xn-nx=0.$$

(1) $\implies$  (2): We obtain $X$ by taking the universal cover $\widetilde{W}$ of $W$ and capping with copies of $D^{4}$. The semi-free action on $X$ is extended from the covering transformations on $\widetilde{W}$. Note that such $X$ is simply-connected.

{(2) $\implies$ (1): We take a metric on $X$ that is invariant under the $\mathbb{Z}/n$ action. By removing geodesic balls surrounding the fixed points, we obtain a free $\mathbb{Z}/n$ action on a 4-manifold $W_0$ with $\partial W_0=\bigcup_{m} S^{3}$. We let $W_{1}$ be the orbit space. Then $\partial W_{1}=\bigcup^{m}_{i=1}L(n,q_{i})$. Moreover, the principal bundle $W_{0}\to W_{1}$ is pulled back from the universal bundle $\widetilde{K}\to K$ via some map $\widetilde{f}: W_{1}\to K$. Here $K$ is the $K(Z_n,1)$-space we constructed and $\widetilde{K}$ is its universal cover. Then we apply Proposition \ref{surgery} to obtain a $\pi_{1}$-isomorphic map $\widetilde{f}':W\to K$ from some 4-manifold $W$ with $\partial W\cong \partial W_{0}$.}  
\end{proof}

\begin{lemma}\label{degree}
Let $D_{iso}(L(n,q), L(n,q'))$ be the set of mapping degrees of those $\pi_1$-isomorphic maps $ f: L(n,q)\to L(n,q')$. Then
$$D_{iso}(L(n,q), L(n,q'))=\{qq'x^2+nk| x,k\in\mathbb{Z}, x \,\,{\rm is \,\,coprime \,\,with\,\,} n\}. $$
\end{lemma}

\begin{proof} We first prove two facts:

(1) $q\in D_{iso}(L(n,q), L(n,1))$ and $q\in D_{iso}(L(n,1), L(n,q))$, and $qq'\in D_{iso}(L(n,q), L(n,q'))$

%(2) every degree $d\in  D_{\pi_1-iso}(L(n,q), L(n,1))$ is coprime to $n$.

(2) each element $d\in  D_{iso}(L(n,q), L(n,q'))$ satisfies $qq'd$ is coprime to $n$ and $qq'd$ is a quadratic residue mod $n$.

To prove (1), recall $L(n,q)=S^3/\tau_{n,q}$ and $L(n,1)=S^3/{\tau_{n,1}}$. The degree $q$ map $\widetilde f_q: S^3\rightarrow S^3$ 
given by $(z,w)\mapsto (z, w^q)$ maps a $\tau_{n,1}$-orbit to a  $\tau_{n,q}$-orbit,
so it  descends to a map $f_q: \textcolor{black}{L(n,1)\rightarrow L(n,q)}$ of degree $q$.  Correspondingly, the degree $q$ map $\widetilde g_q: S^3$ to $S^3$
given by  $(z,w)\mapsto (z^q,w)$  descends to a map $g_q:\textcolor{black}{L(n,q)\rightarrow L(n,1)}$ of degree $q$. Then \textcolor{black}{$f_q\circ g_q$} is a self-map of $L(p,q)$ of degree $q^2$, which is coprime to $n$. Then \textcolor{black}{$f_q\circ g_q$} is a $\pi_1$-isomorphic map by a theorem of \cite{HKWZ}. So $f_q$ is a $\pi_1$-isomorphic map. So we have proved $q\in D_{iso}(L(n,q), L(n,1))$ and \textcolor{black}{$q'\in D_{iso}(L(n,1), L(n,q'))$.}
 Since $D_{iso}(L(n,q), L(n,1))\times D_{iso}(L(n,1), L(n,q'))\subset D_{iso}(L(n,q), L(n,q'))$, so $qq' \in D_{iso}(L(n,q), L(n,q'))$.

To prove (2), suppose there is a map $f: L(n,q)\rightarrow L(n,q')$ of degree $d$. Consider the map $g_{qq'}: L(n,q')\rightarrow L(n,q)$ of 
degree $qq'$ above.
Then the composition $ g_{qq'}  \circ f :L(n,q)\rightarrow L(n,q')\rightarrow L(n,q)$ is a $\pi_1$-isomorphic self map of $L(n,q)$. So $\degr(g_{qq'}  \circ f)=dqq'$ is a quadratic residue mod $n$ and it is coprime to $n$ by a theorem of \cite{HKWZ}.

Now we prove the lemma. 
For each $d\in D_{iso}(L(n,q), L(n,q'))$, $dqq'$ is a quadratic residue modulo $n$ and $dqq'$ is coprime to $n$. So $dqq'=x^2+kn$ for 
$x,k\in \mathbb{Z}$ and $x$ is coprime to $n$ because $x^2$ is coprime to $n$.
So
$d(qq')^2=qq'(x^2+nk)=qq'x^2+n(qq'k).$ That is 
$$d(qq')^2\equiv  qq'x^2 \, mod \, n.$$

Find $q^*$ such that $q^*(qq')=1\, mod \, n$. Then $$d=d(qq')^2q^{*2} \equiv qq'x^2q^{*2}=qq'(xq^*)^2\, mod \, n.$$ That is
$$d=qq'(xq^*)^2+nk'$$ for some $k'\in Z$. 
Note $x,q^*$ are coprime to $n$. So $xq^*$ is also coprime to $n$.
So
$$D_{iso}(L(n,q), L(n,q'))\subseteq\{qq'x^2+nk| x,k\in\mathbb{Z}, x \,\, {\rm is \,\,coprime \,\,with\,\,} n\}. $$
Then we prove the converse.
For each $d$ such that $d=qqx^2+nk$ for $n,k\in\mathbb{Z}$ and $x,n$ coprime. By \cite{HKWZ}, there exists a $\pi_1$-isomorphic map $h: L(n,q')\rightarrow L(n,q')$ of degree $x^2$. $h\circ f:L(n,q)\rightarrow L(n,q')\rightarrow L(n,q')$ has degree $qq'x^2$. Let $(h\circ f)\#p_{n,k}:L(n,q)=L(n,q)\# S^3\rightarrow L(n,q)\wedge S^3\rightarrow L(n,1)$ where $p_{n,k}: S^3\rightarrow L(n,1)$ has degree $-nk$. Then $\degr((h\circ f)\# p_{n,k})=qq'x^2+nk=d$. 
So $(g\circ f)\# p_{n,k}:L(n,q)\rightarrow L(n,q')$ is a $\pi_1$-isomorphic map of degree $d$.
So
$$D_{iso}(L(n,q), L(n,q'))\supseteq\{qq'x^2+nk| x,k\in\mathbb{Z}, x \,\, {\rm is \,\,coprime \,\,with\,\,} n\}.$$
\end{proof}

%\begin{corollary}
%Two lens spaces are $\pi_1$-isomorphic corbordism 
%if and only if  there are types of a semi-free $\mathbb{Z}/n$ action  on a simply connected closed 4-manifolds with $2$ fixed point, and  
%if and only if there is a degree one map from one to another.
%\end{corollary}

\begin{lemma}\label{lem: homotopy equiv between Lens}
The following statements are equivalent:

{\rm (1)} There exists an orientation-preserving homotopy equivalence between 
$L(n,q)$ and $L(n,q')$;

{\rm (2)} There is a degree 1 map $L(n,q)\to L(n,q')$, that is $1\in D_{iso}(L(n,q),L(n,q'))$;

{\rm (3)} There exists integers $x_1$, $x_2$ coprime to $n$ such that $qx_1^2-q'x_2^2\equiv 0\mod n$.
\end{lemma}

\begin{proof}
Any orientation preserving homotopy equivalence is $\pi_{1}$-isomorphic and has mapping degree one. So statement (1) implies statement (2). To see the other direction, let $f: L(n,q)\to L(n,q')$ be a degree-one map. Then $f$ is $\pi_{1}$-surjective and hence $\pi_{1}$-isomorphic. Therefore, one can lift $f$ to degree-one map $\widetilde{f}: S^{3}\to S^{3}$
between their universal covers. By the Hopf theorem, $\widetilde{f}$ is an homotopy equivalence. So $\widetilde{f}$ and $f$ both induce isomorphisms on all higher homotopy groups. By the Whitehead theorem, $f$ is an orientation preserving homotopy equivalence. This shows the equivalence between statement (1) and statement (2).  

Suppose $1\in D_{iso}(L(n,q), L(n,q'))$. Then by Lemma \ref{degree}, there exists $x$ such that 
$1\equiv qq'x^2\mod n.$
Clearly $x$ is coprime to $n$. Then 
$$q(q'x)^2-q'\equiv 0\mod n.$$
Since both $q'$ and $x$ are coprime to $n$, so is $q'x$.  Setting $x_1=q'x$ and $x_2=q'$, we get the congruence relation in (3). This shows that statement (2) implies statement (3).

Suppose there exist $x_1$, $x_2$ coprime to $n$ such that $qx_1^2-q'x_2^2\equiv 0\mod n$. Then
$$qq'x_1^2\equiv (q'x_2)^2\mod n.$$
Since both $q'$ and $x_2$ are coprime to $n$, so is $q'x_2$.
So there exist an integer $l$  such that $lq'x_2\equiv 1\mod n$. Then $qq'(x_1l)^2\equiv 1\mod n$. By Lemma \ref{degree}, we get $1\in D_{iso}(L(n,q), L(n,q'))$. \textcolor{black}{This shows that statement (3) implies statement (2).}
\end{proof}

\begin{proof}[Proof of Theorem \ref{appl-2}] Two lens spaces $L(n,q)$ and $L(n,q')$ are $\pi_1$-isomorphic cobordant if and only if $L(n,q)\bigcup L(n,-q')$ $\pi_{1}$-isomorphicly bounds a 4-manifold $W$. By Theorem \ref{appl1}, this happens if and only if there exist $x_{1},x_{2}$ coprime to $n$ such that $qx^{2}_{1}-q'x^2_{2}\equiv 0 \mod n$. By Lemma \ref{lem: homotopy equiv between Lens}, such $x_{1},x_{2}$ exist exactly when there is an orientation homotopy equivalence between $L(n,q)$ and $L(n,q')$.
\end{proof}

\subsection{Almost free actions} 

We restate Theorem \ref{appl3} as

\begin{theorem}\label{appl 3}
For each 3-manifold $Y\ne S^3$ with $\pi_{1}(Y)$ finite, there exists a finite group $G$, a closed simply connected 4-manifold $X$, and an almost free $G$-action with orbit type $Y$. 
Moreover, we can pick $X$ so that the underlying space of $X/G$ is also simply connected.
\end{theorem}

\begin{proof} By Theorem \ref{main}, we know that $Y$ $\pi_1$-injectively bounds a smooth orientable 4-manifold $W$ with $\pi_{1}(W)$ finite. {We take an injection $\pi_{1}(W)\to A_{n}$ with $n\geq 5$ and consider the corresponding map $W\to K(A_{n},1)$. Since the composition
\[
Y\hookrightarrow W\to K(A_{n},1) 
\]
is $\pi_{1}$-injective and sends $[Y]$ to $0$. We may apply Proposition \ref{surgery} and obtain another manifold $W'$ which is $\pi_{1}$-injectively bounded by $Y$ and has $\pi_{1}(W')=A_{n}$, where $A_n$ is the alternative group of $n$ elements for some large $n$. By replacing $W$ with $W'$, we may assume $\pi_{1}(W)$ is a finite simple group.}

Let $\widetilde W$ be the universal cover of $W$. Then $p: \widetilde W\to W$ is a finite covering with deck transformation
group $G$. Since the inclusion $Y\to W$ is $\pi_1$-injective, it follows that each component $\widetilde Y$ of $p^{-1}(Y)$ is a universal cover of $Y$. So we get $$\partial \widetilde W=p^{-1}(Y)=\sum _{i=1}^n S^3_i,$$
where each $S^3_i$ is a copy of $S^3$.

{Let $X$ be the 4-manifold obtained from $\widetilde{W}$ by capping each boundary component with a 4-ball $B_i^4$. The deck transformation
group $G$ acts freely on $\widetilde W$ with $\widetilde W/G=W$.  Let $G_i\subset G$ be the stabilizer of $S^{3}_i$. Then $G_{i}$ acts on $S^{3}_{i}$ as a covering transformation. So $G_{i}$ is conjugate to a linear action. As a result, this $G$ action can be extended smoothly to $X$.  
Clearly $X$ is simply connected, and 
the $G$ action is almost free with the orbit type $Y$.}

Then  $X/G= CY\cup_{Y}W$, where $CY$ is the cone of $Y$.
Since $CY$ is simply connected and the inclusion $Y\to W$ is $\pi_1$-injective, by Van Kampen Theorem, $\pi_1(X/G)=\pi_1(W)/N$, where $N\subset \pi_1(W)$ is the normal subgroup generated by $\pi_1(M)$.
Since $\pi_1(W)$ is simple and $N$ is non-trivial, we have $\pi_1(W)=N$. Thus $\pi_1(X/G)=1$.
\end{proof}
%\newpage
\section{Minimal bounding index $O_b(Y)$}
\subsection{Finite index bounding and  virtual achirality of aspherical 3-manifolds}

Now we state a more comprehensive version of Theorem \ref{thm: O(Y)=2 or finite}.

\begin{theorem}\label{O(Y)=2 or finite} Let $Y$ be a closed, orientable 3-manifold.
\begin{enumerate}
 \item If $Y$ is aspherical, then $O_{b}(Y)<\infty$ implies that  $Y$ is virtually achiral.
    \item  If $Y$ admits an orientation reversing free involution, then $O_{b}(Y)=2$. The reverse is also true if $Y$ is hyperbolic. 
    
 \item    Suppose $Y$ $\pi_1$-injectively bounds a compact orientable 4-manifold $W$. Then for any integer $d>0$, $Y$ $\pi_1$-injectively bounds a compact orientable 4-manifold $W_d$ such that \[|\pi_1(W_d): \pi_1(Y)|=d|\pi_1(W): \pi_1(Y)|.\]
 \end{enumerate}  
\end{theorem}

We start with some technical lemmas. 

\begin{lemma}\label{lem: Haken cover}
Let $Y$ be an aspherical 3-manifold. Suppose $G$ contains $\pi_{1}(Y)$ as a finite-index subgroup. Then there exists a finite-sheeted covering map $p:\widetilde{Y}\to Y$ such that $\widetilde{Y}$ is Haken and $p_{*}(\pi_{1}(\widetilde{Y}))$ is a normal subgroup of $G$.
\end{lemma}
\begin{proof}
 Since $Y$ is aspherical, there exists a finite cover $p_{1}:Y_{1}\to Y$ such that $Y_{1}$ is Haken \cite{Ag}. Let \[H_{1}=p_{1,*}(\pi_{1}(Y_{1}))\subset \pi_{1}(Y)\subset G.\] Consider the group
\[
H_{2}:=\bigcap\limits_{\gamma\in G} \gamma\cdot H_{1}\cdot \gamma^{-1}.
\]
Then $H_{2}$ is a finite-index normal subgroup of both $H_{1}$ and $G$. Let $p_{2}: Y_{2} \to Y_{1}$ be the covering space that corresponds to  $H_{2}$.  Then $Y_{2}$ is also Haken. The proof is finished by setting $\widetilde{Y}=Y_{2}$ and $p=p_{1}\circ p_{2}$.   \end{proof}

\begin{lemma}\label{lem: orientation reversing deck transformation}
Let $i:J\to G$ be the inclusion of a finite-index normal subgroup. Suppose that $H_{3}(J)=\mathbb{Z}$ and that the map 
\[
i_{*} : H_{3}(J)\to H_{3}(G)
\]
is trivial. Then there exists $\gamma\in G\setminus J$ such that $\phi_{\gamma,*}=-\operatorname{Id}$. Here  
\[\phi_{\gamma,*}: H_{3}(J)\to H_{3}(J)\] is the map induced by automorphism \[\phi_{\gamma}:J\to J,\quad g\mapsto \gamma g\gamma^{-1}.\]
\end{lemma}
\begin{proof} We let $K=K(G,1)$ and let $p:\widetilde{K}\to K$ be the normal covering space that corresponds to $J$. Then \[p_{*}: \mathbb{Z}\cong H_{3}(\widetilde{K})\to H_{3}(K)\]
is trivial because it equals  $i_{*}$. Since $H_{3}(J)=\mathbb{Z}$, each $\phi_{\gamma,*}=\pm \operatorname{Id}$. Suppose $\phi_{\gamma,*}\neq -\operatorname{Id}$ for all $\gamma$. Then the group of deck transformations on $\widetilde{K}$ acts trivially on $H_{3}(\widetilde{K})$. By the universal coefficient theorem, the group of deck transformations also acts trivially on $H^{3}(\widetilde{K};\mathbb{Q})$. Therefore, 
\[p^{*}: H^{3}(K;\mathbb{Q})\to H^{3}(\widetilde{K};\mathbb{Q})\]
is an injection (see \cite[Proposition 3G1]{Ha1}). This is a contradiction.
\end{proof}

\begin{lemma}\label{lem: embedding to Iso(H)}
Let $J$ be  a discrete subgroup of $\operatorname{Iso}(H^3)$ with finite covolume. Suppose $J$ is an index-2 subgroup of $G$ and the inclusion $J\to G$ has no left inverse. Then there exists an embedding $\psi:G\to \operatorname{Iso}(H^3)$ such that $\psi$ sends every element of $J$ to itself and $\psi(G)$ is a finite covolume discrete subgroup of  $\operatorname{Iso}(H^3)$.\end{lemma}
\begin{proof}
Take any $g\in G\setminus J$. Consider the automorphism 
\[
\phi_{g}: J\to J\quad h\mapsto g h g^{-1}.
\]
By Mostow Rigidity, there is a $\gamma \in \operatorname{Iso}(H^3)$ such that $\phi_{g}(h)=\gamma h \gamma^{-1}$ for all $h\in J$. Then we have 
$$g^2hg^{-2}=\phi_{g}^2(h)=\gamma^2h\gamma^{-2}$$ 
for any $h\in J$. That means $\gamma^{2}g^{-2}\in \operatorname{Iso}(H^3)$ commutes with all elements in $J$. Since $J$ is a non-elementary Kleinian group, which implies that 
$\gamma^2=g^2$ \cite[Lemma 1.2.4]{MR}. Consider the 
coset decomposition 
$G=J\sqcup gJ$. We define a map 
$\psi: G\to \operatorname{Iso}(H^3)$
\[
\psi(h)=\begin{cases} h &\text{ if } h\in J,\\
\gamma g^{-1}h &\text{ if } h\notin J.
\end{cases}
\]
Then $\psi$ is a group homomorphism. Since $\psi$ restricts to the identity map on $J$, we have the commutative diagram 
\[
\xymatrix{
J\ar@{=}[d]\ar@{^{(}->}[r] & G\ar@{->>}[d]^{\psi}\ar@{->>}[r] & G/J\cong \mathbb{Z}/2\ar@{->>}[d]^{\psi/J}\\
\psi(J)\ar@{^{(}->}[r]& \psi(G)\ar@{->>}[r] &\psi(G)/\psi(J).
}
\]
Since the inclusion $J\hookrightarrow G$ has no left inverse, $\psi(G)$ must be strictly larger then $\psi(J)$. Hence the group $\psi(G)/\psi(J)$ is nontrivial and the surjective map $\psi/J$ must also be injective. This implies $\psi$ is injective as well.
\end{proof}

%\subsection{ Finite index bounding and achirality}

%\begin{theorem}\label{thm: O(Y)=2 or finite} Let $Y$ be a closed orientable 3-manifold.
%\begin{enumerate}
%\item Suppose $Y$ is aspherical, then $O_{b}(Y)<\infty$ only if $Y$ virtually achiral.
 % \item  If $Y$ admits an orientation reversing free involution, then $O_{b}(Y)=2$. The reverse is also true if $Y$ is hyperbolic. 
 %\end{enumerate}

%\end{theorem}
Now we are ready to prove Theorem \ref{thm: O(Y)=2 or finite}.
\begin{proof}[Proof of Theorem \ref{thm: O(Y)=2 or finite}]
(1) Suppose $W$ is a compact orientable 4-manifold with $\partial W=Y$ and  the inclusion 
$i: Y\to W$ is $\pi_1$-injective and satisfies \[|\pi_1(W): i_*(\pi_1(Y))|<\infty.\] 
Let $G=\pi_{1}(W)$. By Lemma \ref{lem: Haken cover}, there exists a finite cover $p:\widetilde{Y}\to Y$ such that $i_{*}\circ p_{*}(\pi_{1}(Y))$ is a finite index normal subgroup of $G$, denoted by $J$. Let $X$ be a $K(G,1)$-space obtained by attaching cells to $W$. Let $p_{X}:\widetilde{X}\to X$ be the normal covering that corresponds to the subgroup $J$.
Let $f:Y\to X$ be the composition of $Y\xrightarrow{i}W\hookrightarrow X$ and let $\widetilde{f}: \widetilde{Y}\to \widetilde{X}$ be its lift. Then $\widetilde{f}$ is a homotopy equivalence because induces isomorphism between the fundamental group of two apsherical spaces. As a result, we have the following commutative diagram  
\[\xymatrix{
H_3(\widetilde{Y})\ar[rr]^{\widetilde{f}_{*}}_{\cong}\ar[d]^{p_{*}}& & H_{3}(\widetilde{X})\ar[d]^{p_{X,*}}\\
H_3(Y)\ar[r]^{i_{*}=0}&H_3(W)\ar[r] & H_{3}(X).}\]
From this, we see that $p_{X,*}=0$. In other words, the inclusion $J\hookrightarrow G$ induces a trivial map on $H_{3}(-)$. By Lemma \ref{lem: orientation reversing deck transformation}, there exists $\gamma\in G$ such that the automorphsim
$\phi_{\gamma}: J\to J$
induces $-\operatorname{Id}$ on $H_{3}(J;\mathbb{Z})$. Since $J=\pi_{1}(\widetilde{Y})$ and $\widetilde{Y}$ is aspherical, $\phi_{\gamma}$ induces an orientation reversing homotopy equivalence $\tau: \widetilde{Y}\to \widetilde{Y}$. Since $\widetilde{Y}$ is Haken, $\tau$ is homotopic to an orientation reversing homeomorphism. So $Y$ is virtually achiral.

(2)  Suppose $Y$ admits an orientation reversing free involution $\tau$. Then $Y/\tau$ is a closed, non-orientable 3-manifold. Let $W$ be the twisted  $I$-bundle
over $Y/\tau$ {associated to the double cover $Y\to Y/\tau$}. Then $Y$ is the boundary of $W$. The inclusion $Y\to W$ is $\pi_1$-injective and of index $2$. 

Now suppose $Y$ is a hyperbolic  $3$-manifold that $\pi_{1}$-injectively bounds a 4-manifold $W$ with $[\pi_{1}(W): \pi_{1}(Y)]=2$. Let $G=\pi_1(W)$ and let $J=\pi_{1}(Y)$. Then by Proposition \ref{pro: key proposition}, the inclusion $i: J\to G$ induces a trivial map on $H_{3}(-;\mathbb{Z})$. So the inclusion $J\to G$ admits no left inverse. By Lemma \ref{lem: embedding to Iso(H)}, we can regard $G$ a cofinite volume subgroup of $\operatorname{Iso}(H)$ and identify $Y$ with $H^{3}/J$. By Lemma \ref{lem: orientation reversing deck transformation}, there exists some $\gamma\in G\setminus J$ such that the map 
\[
\phi_{\gamma}: J\to J,\quad g\mapsto \gamma g\gamma^{-1}
\]
induces $-\operatorname{Id}$ on $H_{3}(J)=H_{3}(Y)$. In other words, the involution $\tau:Y\to Y$ defined by
\[
\quad [x]\mapsto [\gamma(x)], \quad\forall x\in H^{3}
\]
is orientation reversing. 

It remains to prove $\tau$ is free. Suppose this is not the case. Let $\operatorname{Fix}(\tau)$ be the fixed point of $\tau$. Then we have a decomposition  $\operatorname{Fix}(\tau)=F_0\cup F_2^+ \cup F_2^-$,
where $F_0$ is a union of isolated points, $F_2^+$ and $F_2^-$ are  closed surfaces,  orientable and non-orientable respectively. 

Now we explicitly construct %a space 
a $K(G,1)$-space $P$ and a map $f:Y\to P$ that induces the map $i: J\to G$. Consider $$U=(Y\times [-1, 1])/\widetilde \tau,$$ where $\widetilde \tau(x, t)=(\tau (x), -t)$. Then $U$ is an orbifold with singular loci
\[\operatorname{Fix} (\widetilde \tau)=(F_0\cup F_2^+ \cup F_2^-)\times \{0\}.\] Let $V=U\setminus N$, where
$N$ is an open tubular neighborhood of $\operatorname{Fix} (\widetilde \tau)$. Then $V$ is a manifold with boundary. Other then $Y$, components of $\partial V$ one-to-one correspond to components of $\operatorname{Fix} (\widetilde \tau)$. The space $P$ is obtained by attaching CW complexes to  these components:
\begin{itemize}
    \item Each point in $F_0$ gives a $\mathbb{RP}^{3}$ component of $\partial V$. We attach a copy of $\mathbb{RP}^{\infty}$ via the inclusion   $\mathbb{RP}^{3}\to  \mathbb{RP}^{\infty}$.
    \item Each component of $F^{+}_{2}$ is a closed, orientable surface $F$. It corresponds to a component of $\partial V$ homeomorphic to $F\times \mathbb{RP}^{1}$. We attach a copy of $F\times \mathbb{RP}^{\infty}$ via the standard inclusion $S^{1}=\mathbb{RP}^{1}\to \mathbb{RP}^{\infty}$.
    \item For each component $F'$ of $F^{-}_{2}$, the corresponding component of $\partial V$ is homeomorphic to the unique $\mathbb{RP}^{1}$-bundle over $F'$ whose total space is orientable. We denote it by $F'\widetilde{\times}\mathbb{RP}^{1}$. Since the reflection on $\mathbb{RP}^{1}=S^{1}$ can be extends to an involution of $\mathbb{RP}^{\infty}$, we can define a bundle $F'\widetilde{\times}\mathbb{RP}^{\infty}$ that contains $F'\widetilde{\times}\mathbb{RP}^{\infty}$ as a subbundle. Then we attach a copy of $F'\widetilde{\times}\mathbb{RP}^{\infty}$ to $V$ via the inclusion $F'\widetilde{\times}\mathbb{RP}^{1}\to F'\widetilde{\times}\mathbb{RP}^{\infty}$.
\end{itemize}
Note that $U$ is the quotient of $H^{3}\times [-1,1]$ under a $G$-action. Let $\widetilde{N}$ be the preimage of $N$ under the quotient map $q:H^{3}\times [-1,1]\to U$. Each component of $\widetilde{N}$ is homeomorphic to an open disk, so is contractible. Let $\widetilde{P}$ be the universal cover of $P$. Then $\widetilde{P}$ is obtained by gluing to $(H^{3}\times I)\setminus \widetilde{N}$ 
copies of universal covers of $\mathbb{RP}^{\infty}, F\times \mathbb{RP}^{\infty}, F'\widetilde{\times} \mathbb{RP}^{\infty}$. 
In other words, $\widetilde P$ is obtained by removing contractible subspaces from $H^3\times [-1,1]$ and regluing new contractible spaces. So $\widetilde P$ is homotopy equivalent to $H^3\times [-1,1]$ and $P=\widetilde{P}/G$ is a $K(G,1)$-space.

 For  prime $p$, we use $\mathbb F_p$
for the field of $p$ elements, $\mathbb F_p^\times$ be the  its invertible elements.

\begin{lemma}\label{RP3} Each of the following inclusion map 

{\rm (a)} $\mathbb{RP}^3\rightarrow \mathbb{RP}^\infty$,

{\rm (b)} $F\times \mathbb{RP}^{1}\rightarrow F\times \mathbb{RP}^{\infty}$, 

{\rm (c)} $F'\widetilde\times \mathbb{RP}^{1}\rightarrow F'\widetilde\times \mathbb{RP}^{\infty}$

induces an injection on $H_3(-;\mathbb{F}_2)$.
\end{lemma}

\begin{proof}  (a) is well known. (b) follows from the K\"unneth formula. To prove (c), we consider the Serre spectral sequences for $H_{3}(-;\mathbb{F}_2)$. The only automorphism on $H_*(\mathbb{R}\mathbb{P}^{\infty};\mathbb{F}_2)$ is the identity. So the local coefficients are trivial. For $F'\widetilde\times \mathbb{RP}^{1}$, the differential  $d^2:E^{2}_{2,1}\rightarrow E^2_{0,2}=0$ is trivial. By naturality, the differential 
$d^2:E^2_{2,1}\rightarrow E^2_{0,2}$ for $F'\widetilde\times \mathbb{RP}^{\infty}$ is also trivial. This implies that the map 
$$H_3(F'\widetilde\times \mathbb{RP}^{1};\mathbb{F}_2)\rightarrow H_3(F'\widetilde\times \mathbb{RP}^{\infty};\mathbb{F}_2)$$
is injective.
\end{proof}
Consider the maps on $H_{3}(-;\mathbb{F}_2)$ induced by the inclusions $Y\to V$, $V\to P$ and $Y\to P$. By Lemma \ref{RP3} and the Mayer-Vietoris sequence,  the map \[H_3(V;\mathbb{F}_2)\rightarrow H_3(P;\mathbb{F}_2)\]  is injective. And it is straightforward to see that the map \[H_3(Y;\mathbb{F}_2)\rightarrow H_3(V; \mathbb{F}_2)\] is also injective. So the map 
\[
H_3(Y;\mathbb{F}_2)\rightarrow H_3(P;\mathbb{F}_2)
\]
is injective. However, this is impossible because up to homotopy, the inclusion $Y\to P$ factors through the inclusion $Y\to W$. This contradiction shows that the involution $\tau$ must be free.

(3) Suppose $Y$ $\pi_1$-injectively bounds a compact orientable 4-manifold $W$
with $ [\pi_1(W): \pi_1(Y)]<\infty$. 
Then the inclusion $i:Y\to W$ induce a trivial map on $H_3$.
Given integer $d>1$, we consider the composition
$$f: Y\to W\to W\times L(d,1)\to K(\pi_1(W\times L(d,1)), 1)$$
where the second  map send $W$ to $W\times *$ for some point $*\in L(d,1)$, therefore is $\pi_1$-injective, and the third map is an 
$\pi_1$ isomorphic.
Then $f$ is $\pi_1$-injective, and $f$  induce a trivial map on $H_3$.
Then by Theorem \ref{pro: key proposition}, $f$ $\pi_1$ injectively bounds a compact orientable 4-manifold $W_d$ such that
 $\pi_1(W)\cong \pi_1(W\times L(d,1))$.  Clearly we have \[|\pi_1(W_d): \pi_1(Y)|=d|\pi_1(W): \pi_1(Y)|.\]
 \end{proof}
 
%\subsection{$O_b(L(5,1))$, $\chi_b(L(5,1)$}
%\begin{proposition}\label{lens space} 
%$O_b(L(5,1))=4$ and  $\chi_b(L(5,1))= 2$.
% \end{proposition} 

\subsection{Minimal bounding indices  for lens spaces} 
Let $$d(p)=\min\{d\ge 3|\, d|p-1\}.$$ We restate %the first half of 
Theorem \ref{O(Y)=3} as 

\begin{theorem}\label{O_b(Y)=3}
For each prime $p\ge 5$, $O_b(L(p,q))=d(p).$ 
\end{theorem}

We start with some known facts and  technical lemmas. 

\begin{lemma}\label{22}

{\rm (1)} $H_{2l}(\ZZ_p)=0$,  $H_{2l-1}(\ZZ_p)=\ZZ_p$.

{\rm (2)} $H_{l}(\ZZ_p, \mathbb F_p)=\mathbb F_p$,  $H^{l}(\ZZ_p, \mathbb F_p)=\mathbb F_p$;

{\rm (3)} $H^*( \ZZ_p; \mathbb F_p)\cong\mathbb F_p [x,y]/(x^2)$, with $|x|=1,|y|=2$.
\end{lemma}

\begin{proof} The proof of (1) and (2) are standard calculations in (co)homology of groups.  Calculations of (1) also appear in \cite{SW2}.
(3) is \cite[Chapter XII Section 7]{CE}.
\end{proof}

Recall the universal coefficient theorem 
\begin{equation}\label{eq: UFT1}
0\to H_k(\tilde K)\otimes \mathbb F_p \to H_k(\tilde K, \mathbb F_p) \to {\rm Tor} (H_{k-1}(\tilde K), \mathbb F_p)\to 0    
\end{equation}
and 
\begin{equation}\label{eq: UFT2}
0\to {\rm Ext}(H_{k-1}(\tilde K), \mathbb F_p) \to H_{k}(\tilde K, \mathbb F_p) \to {\rm Hom}(H_{k}(\tilde K), \mathbb F_p)\to 0.    
\end{equation}

Let $\tilde K=K(\ZZ_p, 1)$. Suppose a finite group $D$ acts on $\tilde K$. Then $D$ induces an action $D_k$ on $H_k(\tilde K,\mathbb F_p)=\mathbb F_p$,
which provides representation $\psi_k : D\to \mathbb F_p^\times$, that is, for any $\alpha\in D$, the action of $\alpha$ on $H_k(\tilde K, \mathbb F_p)\cong \mathbb F_p$ is a multiplication by $\psi_k(\alpha)$.

\begin{lemma}\label{24} 
$\psi_{3}(\alpha)=\psi_1(\alpha)^2$ for any $\alpha \in D$.
\end{lemma}

\begin{proof}
By definition, the action of $\alpha$ on $H_1(\tilde K; \mathbb F_p)$ is a multiplication by $\psi_1(\alpha)$. Since 
 $H_0(\tilde K; \mathbb{Z})=\ZZ$, 
by (\ref{eq: UFT1}), there is a natural isomorphism $H_1(\tilde K,\mathbb F_p)\cong H_1(\tilde K, \mathbb{Z})\otimes \mathbb F_p$. %Note that
Since $H_1(\tilde K; \mathbb{Z})\cong \ZZ_p$,
%. Therefore,
 the action of $\alpha$ on $H_1(\tilde K;\mathbb{Z})$ is also a multiplication by $\psi_1(\alpha)$.
Since $H_2(\tilde K, \mathbb{Z})=0$ by (\ref{eq: UFT2}), then by (\ref{eq: UFT2}), we have  $H^2(\tilde K,\mathbb F_p)\cong {\rm Ext}(H_1(\tilde K;\mathbb{Z}), \mathbb F_p)$. Therefore, the action of $\alpha$ on $H^2(\tilde K,\mathbb F_p)$ is also multiplication by $\phi_1(\alpha)$. 
%By Cartan-Eilenberg Theorem, $H^2(\tilde K,F_p)\cong F_p$.
%Let $y$ be a generator of $H^2(\tilde K, F_p)$.
%Note that 

Since $\tilde K=K(\ZZ_p,1)$, we can identify $H_*(\tilde K; \mathbb F_p)$ with $H_*(\ZZ_p,\mathbb{F}_p)$.
By  Lemma \ref{22}, we have
		 $H^2(\tilde K;\mathbb{F}_p)=\langle y\rangle$.
Then the image of $y$ under the action of $\alpha$ equals $\psi_1(\alpha)y$. %by the previous discussion, we can compute that 
Since the action of $\alpha$ preserves cup product, %operation on cohomology,
 the image of $y^2$ under the action of $\alpha$ equals $\psi_1(\alpha)^2y^2$.
Again by Lemma \ref{22}, $H^4(\tilde K,\mathbb F_p)=%F_p
\langle y^2\rangle$. 
Therefore, the action of $\alpha$ on $H^4(\tilde K;\mathbb F_p)$ is multiplication by $\psi_1(\alpha)^2$. 

By Lemma \ref{22}, $H_4(\tilde K;\mathbb{Z})=0$, then by (\ref{eq: UFT2}), we have $H^4(\tilde K,\mathbb F_p)\cong {\rm Ext}(H_3(\tilde K; \mathbb{Z});\mathbb F_p)$.
Since $H_3(\tilde K)=\ZZ_p$, the action of $\alpha$ on $H_3(\tilde K)$ is also a multiplication by $\psi_1(\alpha)^2$.
By Lemma \ref{22} $H_4(\tilde K;\mathbb{Z})=0$, , then by (\ref{eq: UFT1}), we have $H_3(\tilde K,\mathbb F_p)\cong H_3(\tilde K, \mathbb{Z})\otimes \mathbb F_p$. Therefore, the action of $\alpha$ on $H_3(\tilde K; \mathbb F_p)$ is also a multiplication by $\psi_1(\alpha)^2$.
By definition, $\psi_3(\alpha)=\psi_1(\alpha)^2$.
\end{proof}

\begin{lemma}\label{transfer} Let $\pi : \tilde K \to K$ be a finite  regular covering with deck group $D$ and $p$ is a prime.
Suppose $p$ is not a divisor of $|D|$ %and the induced action $D_*$ acts trivially  $H_*(\tilde K, \mathbb F_p)$. 
the induced action of $D$ on $H_*(\tilde K,\mathbb F_p)$ is trivial.
Then $\pi_* : H_*(\tilde K, \mathbb F_p) \to H_*( K, \mathbb F_p)$ 
is an isomorphism. 

Moreover $\pi_* : H_*(\tilde K) \to H_*( K)$ is non-trivial.
\end{lemma}

\begin{proof} In this case, we have the  transfer homomorphism $tr_* : H_*( K, \mathbb F_p)\to H_*(\tilde  K,\mathbb F_p)$ and that the composition
\[\pi_* \circ \tr_* : H_*( K, \mathbb F_p)\to H_*(\tilde  K, \mathbb F_p)  \to H_*( K, \mathbb F_p)\]
is the multiplication by $d=|D|$, that is for each $u\in H_*( K, \mathbb F_p)$, $\pi_* \circ \tr_*(u)=d \cdot u$, for detail, see \cite[p.392]{Ha1}.
Since the induced action of $D$ on $H_*(\tilde K,\mathbb F_p)$ is trivial,
%$D_*$ acts trivially in $H_*(\tilde  K, \mathbb F_p)$, 
it is also easy to verify that 
$ \tr_*\circ \pi_* (v) =d\cdot v$ for each $v\in H_*( \tilde K, \mathbb F_p).$
Since $p$ is not a divisor of $d$,  $d\ne 0$. Let $\bar \tr_*=\tr_*/d$, then 
\[ \bar \tr_*\circ \pi_* =\id, \ \ \pi_* \circ \bar \tr_*=\id, \]
that is $\pi_*: H_*(\tilde K, \mathbb F_p) \to H_*( K,\mathbb F_p)$ is an isomorphism. 

The "Moreover" part:  By  (\ref{eq: UFT1}), we have 
\[\xymatrix{
0\ar[r]&H_3(\tilde K)\otimes \mathbb F_p \ar[r] \ar[d]^{\pi_*} &H_3(\tilde K, \mathbb F_p) \ar[r] \ar[d]^{\pi_*}&{\rm Tor} (H_2(\tilde K), \mathbb F_p)\ar[r]\ar[d]^{\pi_*} &0\\
0\ar[r]&H_3( K)\otimes \mathbb F_p\ar[r] &H_3( K, \mathbb F_p)\ar[r]&{\rm Tor} (H_2( K), \mathbb F_p)\ar[r] &0.
}\]
If $\pi_* : H_*(\tilde K) \to H_*( K)$ is trivial, then  $\pi_* : H_*(\tilde K)\otimes \mathbb F_p  \to H_*( K)\otimes \mathbb F_p$ is trivial.
Since $H_2(\tilde K)=0$, which will contradicts that $\pi_* : H_*(\tilde K, \mathbb F_p) \to H_*( K, \mathbb F_p)$ is an isomorphism. 
\end{proof}

\begin{lemma}\label{11}
There is a group $G$ of order $pd(p)$ and an injection $i: \ZZ_p\to G$
such that the induced map $i_* : H_3(\ZZ_p)\to H_3(G)$ is trivial.
\end{lemma}

\begin{proof} Since $d(p)|p-1$, we can take  $u\in\mathbb F_p^\times$ such that $\text{ord}(u)=d(p)$. Let \[G=\langle\alpha, \beta | \alpha^p=\beta^{d(p)}=1, \, \beta\alpha \beta^{-1} = \alpha^u \rangle.\]
Then $G=\ZZ_p\langle \alpha \rangle \rtimes \ZZ_{d(p)}\langle \beta \rangle$ is a group of order $pd(p)$. Let \[c(\beta): G\to G\ \ \text{be given by} \, x\to \beta x\beta^{-1}.\]
Then $c(\beta)$ keeps $\ZZ_p$ invariant and its restriction on $\ZZ_p$  is $m: \ZZ_p\to \ZZ_p$ given by $\alpha\mapsto \alpha^u$. As an inner automorphism on $G$, $c(\beta)_*$ induces the identity on $H_*(G)$. 
%For the commutative diagram
%$$\begin{CD}
%			\ZZ_p @>  i  >> G   \\
%		         @VV m V   @VV c(\beta) V   \\
%			\ZZ_p @ > i >>  G.
%		\end{CD} $$
Note that $m_*:H_1(Z_p)\to H_1(Z_p)$ is a multiplication by $u$.
By a similar argument as in Lemma \ref{24}, we have $m_*:H_3(Z_p)\to H_3(Z_p)$ is a multiplication by $u^2$. Consider the following diagram on $H_3$:
%Take $H_1$, we have 
%$$\begin{CD}
%			\ZZ_p=H_1(\ZZ_p) @>  i_*  >> H_1(G)   \\
%		         @VV m_* V   @VV \text{id}_* V   \\
%			\ZZ_p=H_1(\ZZ_p) @ > i_* >>  H_1(G).
%		\end{CD} $$
%By Lemma \ref{22},  $H_2(\ZZ_p)=0$. Then  by the Universal Coefficient Theorem, we have
%$$\begin{CD}
%			\ZZ_p=H^2(\ZZ_p,\mathbb F_p) @>  i_*  >> H^2(G,\mathbb F_p)   \\
%		         @VV m_* V   @VV \text{Id}_* V   \\
%			\ZZ_p=H^2(\ZZ_p, \mathbb F_p) @ > i_* >>  H^2(G, \mathbb F_p).
%		\end{CD} $$
%By Lemma \ref{22} (3), we have
%$$\begin{CD}
%			\ZZ_p=H^4(\ZZ_p, \mathbb F_p) @>  i_*  >> H^4(G, \mathbb F_p)   \\
%		         @VV m_* V   @VV \text{Id}_* V   \\
%			\ZZ_p=H^4(\ZZ_p, \mathbb F_p) @ > i_* >>  H^4(G, \mathbb F_p).
%		\end{CD} $$
%and $m(z)= u^2z$, where $z=y^2$ is a generator of  $H^4(\ZZ_p, \mathbb F_p)$. By Lemma \ref{1} (1) $H_4(\ZZ_p)=0$, and then by the Universal Coefficient Theorem, we have
	$$\begin{CD}
			\ZZ_p=H_3(\ZZ_p) @>  i_*  >> H_3(G)   \\
		         @VV m_* V   @VV c(\beta)_*=\text{Id} V   \\
			\ZZ_p=H_3(\ZZ_p) @ > i_* >>  H_3(G).
		\end{CD} $$
and 	$m(w)= u^2w$, where $w$ is a generator of  $H_3(\ZZ_p)$.		So we have 
$$u^2i_*(w)=i_*(u^2w)=i_* m_* (w)=i_*(w).$$
Since $\text{ord}(u)=d(p)\ge 3$, we have $u\ne \pm 1$, so $u^2\ne 1$, and we conclude that $i_*(w)=0$, that is $i_*$ is trivial.
\end{proof}

\begin{proof}[Proof of Theorem \ref{O_b(Y)=3}] We first prove that $O_b(L(p,q))\le d(p)$:
Consider the $\pi_1$-injective map
$$h=f_\psi \circ i: L(p,q)\to K(\ZZ_p,1 ) \to K(G, 1),$$ where $ i: L(p,q)\to K(\ZZ_p,1 ) $ is the inclusion, and $f_\psi$ realizes 
the injection $\psi: \ZZ_p\to G$  on $\pi_1$
given by Lemma \ref{11}. Then 
$$h_*={f_\psi}_* \circ i_*: H_3(L(p,q))\to H_3(K(\ZZ_p,1)) \to H_3(K(G, 1))$$
is a trivial map by Lemma \ref{11}. 
Then by Theorem  \ref{pro: key proposition}, there exists a smooth 4-manifold $W$ bounded by $L(p,q)$, and an isomorphism $\pi_{1}(W)\cong G$ under which $\psi$ is exactly the map induced by the inclusion $L(p,q)\to W$. Now $|\pi_1(W): \ZZ_p|=d(p)$.
Hence  $O_b(L(p,q))\le d(p)$.

Next  we prove  $O_b(L(p,q))\ge d(p)$. Otherwise there is compact 4-manifold $W$ such that $\partial W=L(p,q)$, the inclusion $i: L(p,q)\to W$ is 
$\pi_1$-injective and $|\pi_1(W): \ZZ_p|< d(p)$. By Sylow Theorem, $\ZZ_p$ is a normal subgroup of $G=\pi_1(W)$.
Then we have the regular covering $\pi: \tilde K\to K=K(G, 1)$ with deck group $D=G/\ZZ_p$ and the following commutative diagram
up to homotopy

$$\begin{CD}
			L(p,q)=\partial W @>  i  >> W   \\
		         @VV j V   @VV j' V   \\
			K(\ZZ_p, 1)=\tilde K@> \pi >>  K=K(G, 1).
		\end{CD} $$

 Then we have commutative diagram
 
$$\begin{CD}
			H_3(L(p,q)) @>  i_*  >> H_3(W)   \\
		         @VV j_* V   @VV j'_* V   \\
			H_3(\tilde K) @> \pi_* >>  H_3(K).
		\end{CD} $$
		
		Clearly $i_*$ is a trivial map. 
		Since $j_*$ is a surjection, $\pi_*$ is a trivial map. 

On the other hand $D$ induces an action $D_k$ on $H_k(\tilde K,\mathbb F_p)=\mathbb F_p$,
therefore provides representation $\psi_k : D\to \mathbb F_p^\times$,  
which implies that $|{\rm Im} \psi_1|$ is a divisor of both $|D|$ and $p-1$, therefore a divisor of ${\rm gcd}(|D|, p-1)$.
Since $|D|< d(p)$, it follows that ${\rm gcd}(|D|, p-1)\le 2$, that is $\psi_1(\alpha)=\pm 1$ for any $\alpha \in D$.
By Lemma \ref{24}, $\psi_{3}(\alpha)=\psi_1(\alpha)^2=1$ for any $\alpha \in D$, that is $D$ acts trivially on $H_3(\tilde K, \mathbb F_p)$.
Then $\pi_* : H_*(\tilde K, \mathbb F_p) \to H_*( K,\mathbb F_p)$ 
is an isomorphism and  $\pi_* : H_*(\tilde K) \to H_*( K)$ is nontrivial by Lemma \ref{transfer}.
We reach a contradiction.		
\end{proof}

%The second half of Theorem \ref{O(Y)=3} is contained in the following

%\begin{corollary}
%Suppose $q$ is $1$ or a prime. Then $q$ is realized as $O_b(Y)$ for infintely many $Y$.
%\end{corollary}
%\begin{proof}
%For prime $q=1$ or  2, the proof is contained in Remark \ref{finite index}.
%Now suppose $q$ is prime $\ge 3$. Let $q_1< q_2<... <q_k$ be all  primes $<q$. Then by Chinese Remainder Theorem and by
%Dirichlet Theorem in number theory, the %following congruence equation system 
%\[ p\, \equiv \, -1 \, \text{mod} \, q_1q_2...q_k,\,  \  \ p\, \equiv \, 1\, \text{mod} \,q, \]
%has infinitely many  prime solutions $p$, and then for each such prime $p$, $O_b(L(p, 1))=d(p)=q$ by Theorem \ref{O_b(Y)=3}.
%\end{proof}

\section{Some explicit examples}

\subsection{On surface bundles  bounding  surface bundles} We prove  Proposition \ref{bundle} and Corollary \ref{AW} in this subsbection, and we restate them:

\begin{proposition}\label{bundle1}  Suppose $Y$ is a $\Sigma_g$-bundle over $S^1$, $g\ge 3$. 
Then $Y$ bounds a surface bundle over surfece.
Moreover, the bounding  is $\pi_1$-injective and $W$ has residually finite $\pi_1$.
\end{proposition}

\begin{proof} Let  $\rm{MCG}_+(\Sigma_g)$ be the oriented mapping class 
group of $\Sigma_g$, and $\rm{MCG}_+(\Sigma_g)^\text{ab}$ be its abelianization. 
  Each $\Sigma_g$-bundle over $S^1$ has the form $(\Sigma_g, h)$, where $h: \Sigma_g \to \Sigma_g$  is a homeomorphism. 
Let $Y$ be such a $\Sigma_g$-bundle $(\Sigma_g, h)$. Then $Y$ is pulled back from the universal surface bundle $\Sigma_{g}\hookrightarrow E\to \operatorname{BHomeo}_{+}(\Sigma_{g})$
via a map 
\[
f: S^{1}\to \operatorname{BHomeo}_{+}(\Sigma_{g}).
\]
Here $\operatorname{BHomeo}_{+}(\Sigma_{g})$ is the classifying space of the group of orientation preserving homeomorphisms on $\Sigma_{g}$, which is a $K(\MCG_{+}(\Sigma_{g}), 1)$ space \cite[Section 5.6]{FM}. 
Note that we have 
\[
H_{1}(\operatorname{BHomeo}_{+}(\Sigma_{g}))\cong \pi_{1}(\operatorname{BHomeo}_{+}(\Sigma_{g}))^{\operatorname{ab}}=\MCG_{+}(\Sigma_{g})^{\operatorname{ab}}.
\]
As proved by Mumford \cite{Mum} and Powell \cite{Po}, $\MCG_{+}(\Sigma_{g})^{\operatorname{ab}}=0$ for $g\ge 3$.
 Then $[h]=0\in \MCG_{+}(\Sigma_{g})^{\operatorname{ab}}$, which  implies that $[f]=0\in H_{1}(\operatorname{BHomeo}_{+}(\Sigma_{g}))$. This implies that the map 
\[
f:  S^{1}\to \operatorname{BHomeo}_{+}(\Sigma_{g})
\]
can be extended to a map 
\[
\widetilde{f}: \Sigma_{n,1}\to \operatorname{BHomeo}_{+}(\Sigma_{g}).
\]
Here $\Sigma_{n,1}$ is a surface of genus $n$ with $1$ boundary components. We may assume $n>0$ since otherwise we can precompose $\widetilde{f}$ with a degree-1 map $\Sigma_{1,1}\to \Sigma_{0,1}$. Pulling back the universal bundle over $\operatorname{BHomeo}_{+}(\Sigma_{g})$ via $\widetilde{f}$, we obtain a surface bundle 
\[
\Sigma_{g}\hookrightarrow W\to \Sigma_{n,1}.
\]
The boundary of the total space $W$ is $Y$.

Since $g>0$, $n>0$,  we have  the  commutative diagram 
\[\xymatrix{
1\ar[r]&\pi_{1}(\Sigma_{g})\ar[r] \ar@{=}[d]&\pi_{1}(Y)\ar[r]\ar[d]^{i_{*}}&\pi_{1}(S^{1})\ar[r]\ar@{^{(}->}[d] &1\\
1\ar[r]&\pi_{1}(\Sigma_{g})\ar[r] &\pi_{1}(W)\ar[r]&\pi_{1}(\Sigma_{n,1})\ar[r] &1
}\]
which directly implies that $Y\hookrightarrow W$ is $\pi_{1}$-injective. Note that $\pi_{1}(\Sigma_{n,1})$ is a free group, so the exact sequence 
\[
1\to \pi_{1}(\Sigma_{g})\to \pi_{1}(W)\to \pi_{1}(\Sigma_{n,1})\to 1
\]has a section and implies the isomorphism 
\[
\pi_{1}(W)\cong \pi_{1}(\Sigma_{g})\rtimes \pi_{1}(\Sigma_{n,1}).
\]
Since a semi-direct product of residually finite groups is residually finite, we see that $\pi_{1}(W)$ is residually finite
follows from this diagram. 
Then the  "Moreover" part follows.
\end{proof}

\begin{corollary}\label{AW1}
Suppose $Y$ is a   closed orientable hyperbolic or mixed $3$-manifold. Then  a finite cover of $Y$ bounds a surface bundle over surface.
\end{corollary}
\begin{proof}  By theorems on virtually fibrations of hyperbolic 3-manifolds \cite{Ag} and mixed 3-manifolds \cite{PWi}, 
 $Y$ has a finite cover which is an orientable surface $\Sigma_g$-bundle over a circle with  $g\ge 3$.
 Then Corollary  \ref{AW1} follows from Proposition \ref{bundle1}.
 \end{proof}

\subsection{4-manifolds bounded by $L(5,1)$ realizing $O_b$ and with minimal $\chi$}

%We construct a 4-manifold $\pi_1$-injectively bounded by $L(5, 1)$.
Let $f: \mathbb{CP}^{2}\to \mathbb{CP}^{2}$ be a projective transformation in $PGL_3(\mathbb{C})$ defined as 
$$f([x_1: x_2:  x_3])=[\bar \zeta  x_1: x_2: \zeta x_3],$$
where $\zeta=e^{\frac {2 \pi i}{5}}$. Then $f$ is a generator of $\mathbb{Z}_5$-action on $\mathbb{CP}^{2}$ and has three fixed points
with homogeneous coordinates 
\[P_1=[1: 0: 0],\quad  P_2=[0: 1: 0],\quad P_3=[0: 0: 1].\]
Moreover one can check that they have types $L(5, 2)$,  $L(5, -1)$ and  $L(5, 2)$ respectively. Recall that $L(p,-q)$ is the orientation reversal of $L(p,q)$. Here the type of a fixed point $P$ is defined to be the oriented spherical manifold $\partial D/f$ where $D$ is an $f$-invariant small regular neighborhood of $P$.

%Let $U_1$ be the affine space given by $U_1=\{[x_1:x_2:x_3] \in \mathbb{CP}^{2}\, |\,  x_1\ne 0\}$.
% Let $\phi_1: U\rightarrow \mathbb{C}^2$ given by 
%$\phi_1([x_1:x_2:x_3])=(\frac{x_2}{x_1},\frac{x_3}{x_1})$. Then $(U_1, \phi_1)$ is a chart of $\mathbb{CP}^{2}$. Under this chart, $f$ has representation $\widetilde f_{(U_1,\phi_1)}(x_2,x_3)=\phi_1\circ f\circ \phi_1^{-1}(x_2, x_3)=(\zeta x_2, \zeta^2 x_3). $ So the orbit type at $P_1=[1:0:0]$ is $L(5,2)$.

%Let $U_2$ be the affine space given by $U_2=\{[x_1:x_2:x_3]\in \mathbb{CP}^{2}|x_2\ne 0\}$.
%Let $\phi_2: U_2\rightarrow \mathbb{C}^2$ given by $\phi_2([x_1:x_2:x_3])=(\frac{x_1}{x_2},\frac{x_3}{x_2})$. Then $(U_2,\phi_2)$ is a chart of $\mathbb{CP}^{2}$. Under this chart, $f$ has representation $\widetilde f_{(U_2,\phi_2)}(x_1,x_3)=\phi_2\circ f\circ \phi_2^{-1}(x_1,x_3)=(\zeta^{-1}x_1,\zeta x_3).$
%So the orbit type at $P_2=[0:1:0]$ is $L(5,-1)$.

%Let $\iota: \mathbb{CP}^{2}\rightarrow \mathbb{CP}^{2}$ be the involution given by $\iota([x_1:x_2:x_3])=[x_3:x_2:x_1]$. Then $f^{-1}=\iota^{-1}\circ f\circ \iota$ and $\iota(P_1)=P_3$. So the orbit type at $P_3=[0:0:1]$ is isomorphic to the orbit type at $P_1$. So the orbit type at $P_3$ is also $L(5,2)$.

Let $D_1, D_2, D_3$ be the $f$-invariant regular neighborhoods of $P_1$, $P_2$, $P_3$. Let
$$W=\frac {\mathbb{CP}^{2} \setminus \sum _{i=1}^3 D_i}{f}.$$
Then $\partial W=L(5, -2)\cup L(5, 1)\cup L(5, -2)$ with the induced orientation. Note $\pi_1(W)=\langle\alpha |\alpha^5=1\rangle$, 
the inclusion of each component of $\partial W$ 
into $W$ is $\pi_1$-isomorphic. Gluing two $L(5,2)$ in $\partial W$ via an orientation reversing homeomorphism
(recall $L(5,2)$ admits such homeomorphism \cite{Ha2}), we get a compact oriented  4-manifold $W_1$ bounded by $L(5,1)$ and we can compute the fundamental group of $W_1$ by the HNN-extension theorem:
$$\pi_1(W_1)=\langle\alpha, t| \alpha^5=1, t\alpha t^{-1}=\alpha^r\rangle.$$
Here $1\le r\le 4.$
Let $c$ be a simple closed curve such that the  algebraic intersection number of $c$ and  $L(5,2)\subset W_1$ is 4.
Let $S^1\times D^3$ be a regular neighborhood of $c$.  Doing surgery along $c$, we get a closed oriented 4-manifold
$$W_2=(W_1\setminus S^1\times D^3)\cup D^2\times S^2.$$ 
  
Now we have by Seifert-Van Kampen theorem:
$$\pi_1(W_2)=\langle\alpha, \tau| \alpha^5=1, \tau\alpha \tau^{-1}=\alpha^r,  \tau^4=1\rangle.$$
Consider the automorphism \[\phi: \mathbb{Z}_5\langle\alpha\rangle \rightarrow \mathbb{Z}_5\langle\alpha\rangle,\ \alpha\mapsto\alpha^r.\]
Since $r^4\equiv1({\rm mod}\,\, 5)$ (Fermat's little theorem), the order of $\phi$ divides $4$.
So there is a well-defined homomorphism $\rho: \mathbb{Z}_4\langle\tau\rangle\rightarrow \rm{Aut}(\mathbb{Z}_5\langle\alpha\rangle)$ such that $\rho(\tau)=\phi$.
Then \[\pi_1(W_2)=\mathbb{Z}_5\langle\alpha\rangle\rtimes_\rho \mathbb{Z}_4\langle\tau\rangle\] is a semi-direct product. So $\alpha$ is nontrivial in $\pi_1(W_2)$. So the inclusion map $L(5,1)\rightarrow W_2$ is $\pi_1$-injective.

The order of $\pi_1(W_2)$ is $5\times 4=20$. Since $O_b(L(5,1))=4$ by Theorem \ref{O_b(Y)=3}, $W_2$ realizes $O_b(L(5,1))$.

 We now verify that $\chi_b(L(5,1))=2$. It is easy to see that $\chi(\mathbb{CP}^{2})=3$, so $\chi(\mathbb{CP}^{2} \setminus \sum _{i=1}^3 D_i)=0$, and then $\chi(W)=0$.
 Since $\chi(Y)=0$ for each closed 3-manifold, by the gluing formula of $\chi$, it is easy to see that $\chi(W_1)=0$, and then 
 $\chi(W_2)=2$. So we have $1\le \chi_b(L(5,1))\le 2$. 
 
 Suppose $\chi_b(L(5,1))=1$. Then $L(5,1)$ bounds a rational homology 4-ball $W'$. This is impossible because $|H_1(L(5,1);\mathbb{Z})|=5$ is not a square number (see \cite[Lemma 3]{CG}). So $\chi_b(L(5,1))=2$. Therefore $W_2$ realizes $\chi_b(L(5,1))$.

 % \begin{lemma}\label{appl4}
%For each closed orientable 3-manifold $Y$ with finite $\pi_1$, $\chi_b(Y)\ge 1$. 
%\end{lemma}
%\begin{proof}

%Since $|\pi_1(W)|<\infty$, the first Betti number $b_1(W)=0$. Since $|\pi_1(\partial W)|<\infty$, we have 
%$H^1(W, Q)=0$. Then by the cohomolgy exact sequence of the pair $(W, \partial W)$, we have $H^3(W)\cong H^3(W,Y)$,
%and by Lefschetz duality  $0=H_1(W)\cong H^3(W,Y)$, we have $H^3(W)=0$.
%So $\chi(W)=1+b_2(W)\ge 1$, so  $\chi_b(Y)\ge 1$.
%\end{proof}

%The proof of the main result, Theorem \ref{main}, relies on several ingredients including Thurston's picture on 3-manifolds, and  the virtually fiberation theorems
%of Agol and of Przytycki-Wise
%(Theorem \ref{v fiber}). as well a fact derived from mapping class groups of surfaces (Proposition \ref{fiber torsion}. We start by recalling the following deep theorem. 

\end{document}